\documentclass[11pt,leqno]{article}
\usepackage{amsthm,amsfonts,amssymb,amsmath,oldgerm}
\usepackage{epsfig}
\numberwithin{equation}{section}
\usepackage[thinlines]{easybmat}

\newcommand{\DDD}{D3'}

\newcommand{\bu}{\bar{U}}

\def\eps{\varepsilon }

\newcommand\R{\mathbb R}

\def\eps{\varepsilon}

\newcommand\br{\begin{remark}}
\newcommand\er{\end{remark}}
\newcommand\bp{\begin{pmatrix}}
\newcommand\ep{\end{pmatrix}}
\newcommand\be{\begin{equation}}
\newcommand\ee{\end{equation}}
\newcommand\ba{\begin{equation}\begin{aligned}}
\newcommand\ea{\end{aligned}\end{equation}}

\newcommand{\bap}{\begin{app}}
\newcommand{\eap}{\end{app}}
\newcommand{\begs}{\begin{exams}}
\newcommand{\eegs}{\end{exams}}
\newcommand{\beg}{\begin{example}}
\newcommand{\eeg}{\end{exaplem}}
\newcommand{\bpr}{\begin{proposition}}
\newcommand{\epr}{\end{proposition}}
\newcommand{\bt}{\begin{theorem}}
\newcommand{\et}{\end{theorem}}
\newcommand{\bc}{\begin{corollary}}
\newcommand{\ec}{\end{corollary}}
\newcommand{\bl}{\begin{lemma}}
\newcommand{\el}{\end{lemma}}
\newcommand{\bd}{\begin{definition}}
\newcommand{\ed}{\end{definition}}
\newcommand{\brs}{\begin{remarks}}
\newcommand{\ers}{\end{remarks}}

\newtheorem{theo}{Theorem}[section]

\newtheorem{exams}[theo]{Examples}

\numberwithin{equation}{section}

\newcommand{\RR}{{\mathbb R}}

\newtheorem{theorem}{Theorem}[section]
\newtheorem{proposition}[theorem]{Proposition}
\newtheorem{corollary}[theorem]{Corollary}
\newtheorem{lemma}[theorem]{Lemma}
\newtheorem{definition}[theorem]{Definition}

\newtheorem{example}[theorem]{Example}
\newtheorem{remark}[theorem]{Remark}

\newtheorem{thm}{Theorem}

\newcommand{\RM}{\mathbb{R}}
\newcommand{\ZM}{\mathbb{Z}}

\newcommand{\CM}{\mathbb{C}}

\title
{Whitham averaged equations and modulational stability of periodic traveling waves of a  hyperbolic-parabolic balance law}

\author{\sc \small
Blake Barker\thanks{Indiana University, Bloomington, IN 47405;
bhbarker@indiana.edu: Research of B.B. was partially supported
under NSF grants no. DMS-0300487 and DMS-0801745.}
~~~
Mathew A. Johnson\thanks{Indiana University, Bloomington, IN 47405;
matjohn@indiana.edu: Research of M.J. was partially supported by an NSF Postdoctoral Fellowship under NSF grant DMS-0902192.}
~~~
Pascal Noble\thanks{Universit\'e Lyon I, Villeurbanne, France;
noble@math.univ-lyon1.fr:
Research of P.N. was partially supported by the French ANR Project no.
ANR-09-JCJC-0103-01.}
\\
\sc \small
~~~
L.Miguel Rodrigues\thanks{Universit\'e de Lyon, Universit\'e Lyon 1,
Institut Camille Jordan, UMR CNRS 5208, 43 bd du 11 novembre 1918,
F - 69622 Villeurbanne Cedex, France; rodrigues@math.univ-lyon1.fr}
~~~
Kevin Zumbrun\thanks{Indiana University, Bloomington, IN 47405;
kzumbrun@indiana.edu:
Research of K.Z. was partially supported
under NSF grants no. DMS-0300487 and DMS-0801745.}
}
\begin{document}

\maketitle


\begin{center}
\end{center}

\begin{center}
{\bf 2000 MR Subject Classification}: 35B35.
\end{center}


\begin{abstract}
In this note, we report on recent findings concerning the spectral and nonlinear
stability of periodic traveling wave solutions of
hyperbolic-parabolic systems of balance laws,
as applied to the St. Venant equations of shallow water flow down an
incline.
We begin by introducing a natural set of spectral
stability assumptions, motivated by considerations from the Whitham averaged equations,
and outline the recent proof yielding nonlinear stability under these conditions.  We then
turn to an analytical and numerical investigation of the verification of these spectral stability assumptions.
While spectral instability is shown analytically to hold in both the Hopf and homoclinic limits,
our numerical studies indicates spectrally stable periodic solutions of intermediate period.
A mechanism for this
moderate-amplitude stabilization is proposed in terms of
numerically observed  ``metastability" of the
the limiting homoclinic orbits.
\end{abstract}
\newpage


\bigbreak

\section{Introduction }\label{intro}
Nonclassical viscous conservation or balance laws arise in many areas of mathematical modeling including
the analysis of multiphase fluids or solid mechanics.  Such equations are known to exhibit a wide variety of traveling wave phenomena such
as homoclinic or heteroclinc solutions, corresponding to the standard pulse and front or shock type solutions, respectively,
as well as solutions which are spatially periodic.  Historically, a great deal of effort has been applied to understanding
the time-evolutionary stability of the homoclinic/heteroclinic solutions of such equation, and their spectral and
nonlinear stability theories are well understood (in general).  In contrast, until recently the analogous stability
theories of the periodic counterparts have received relatively little attention.  The goal of this paper is to
present recent progress towards the understanding of periodic
traveling-wave solutions within the context of a particular physically
interesting
hyperbolic--parabolic system of second order PDE's which we describe below.  We begin, however, by briefly
recalling the known theory in the case of a strictly parabolic system
of conservation laws.

There has been a great deal of recent progress towards the understanding of the stability properties of periodic
traveling waves of viscous strictly parabolic systems of conservation laws 
of the form
\begin{equation}\label{e:conslaw-md}
u_t+ \nabla\cdot f(u)=\Delta u, \quad x\in\RM^d, u\in \RM^n;
\end{equation}
see \cite{JZ4,JZ2,OZ1,OZ2,OZ3,OZ4,Se1}.
In particular, using delicate analysis of the resolvent of the linearized operator it has been shown that any periodic
traveling wave solution of \eqref{e:conslaw-md} that is spectrally stable
with respect to localized ($L^2$) perturbations is
(time-evolutionary) nonlinearly stable in $L^p$ or $H^s$ for appropriate values of $p$ and $s$.  In fact, such solutions are
asymptotically stable (in an appropriate sense) for dimensions $d\geq 2$, while they are only nonlinearly bounded stable
for dimensions $d=1$.
%
However, while these results are mathematically satisfying,
up to now no example of a spectrally stable periodic solution of
equations of the form \eqref{e:conslaw-md} has yet been found.
In fact, in dimension one it was shown in \cite{OZ1}
by rigorous Evans function computations that such solutions cannot exist for
certain
model systems of form \eqref{e:conslaw-md} admitting a Hamiltonian structure,
due to the existence of a spectral dichotomy.

%
Though these isolated
results may seem discouraging, it should be noted that
the
explicit examples of form \eqref{e:conslaw-md}
considered so far have possibly had too much
or the wrong sort of
 structure to admit stable periodic solution.  In particular, it may very well
be the case that by considering more exotic potentials\footnote{
In particular, potentials for which the period is decreasing with
respect to amplitude in some regime.}
it may be possible to find a stable periodic solution
even
within
the class of examples studied in \cite{OZ1}.
And, for higher dimensional ($n\ge 3$) or more general
systems,
one may expect to find richer behavior, including, possibly,
stable periodic waves.
Whether or not this
occurs for systems of form \eqref{e:conslaw-md}
remains an interesting open problem.
However, within the slightly wider class of nonstrictly parabolic systems
of balance laws,
it has recently been shown that stable periodic waves {\it can indeed occur},
which leads us to the investigations presented in this paper.

A particular class of examples of more general (non-strictly parabolic) balance laws which have been recently seen (numerically)
to admit stable periodic traveling wave solutions is given by the generalized St. Venant equations, which in Eulerian coordinates
take the form
%
\ba \label{e:stv}
h_t + (hu)_x&= 0,\\
(hu)_t+ (h^2/2F+ hu^2)_x&= h- u|u|^{r-1}/h^s +\nu (hu_x)_x ,
\ea
where $1\leq r\leq 2$ and $0\leq s\leq 2$.  Unlike \eqref{e:conslaw-md}, the St. Venant equations are a second order
hyperbolic--parabolic system of PDE's in balance law form (due to the non-differentiated source term in the parabolic
part).  Equations of this form arise naturally when approximating shallow water flow on an inclined ramp,
in which case $h$ represents the height of the fluid, $u$ the velocity average with respect to height, $\nu$ is a nondimensional
viscosity equal to the inverse of the Reynolds number, and $F$ is the Froude number, which here is the square
of the ratio between the speed of the fluid and the speed of gravity waves\footnote{In particular, the Froude number $F$
depends on the angle of inclination of the incline.}.  Further, the term $u|u|^{r-1}/h^s$ models
turbulent friction along the bottom surface and $x$ measures longitudinal distance along the ramp.  Finally, we point out
that the form of the viscosity term $\nu(hu_x)_x$ is motivated by the formal derivations from the Navier Stokes equations
with free surfaces; other choices are obviously available and are sometimes used in the literature \cite{HC}.
Typical choices for the parameters $(r,s)$ are $r\in\{1,2\}$ and $s\in\{0,1,2\}$; see \cite{BM,N1,N2} and references
therein.  The choice $(r,s)=(2,0)$ is considered in detail in the works \cite{N1,N2,JZN,BJRZ}.

Periodic and solitary traveling waves, known as
roll waves,
are well-known to appear as solutions of \eqref{e:stv},
generated by competition between gravitational force and
friction along the bottom.
Such patterns have been used to model phenomena in
several areas of the engineering literature, including landslides, river and spillway flow, and the topography of
sand dunes and sea beds and their stability properties have been much studied
numerically, experimentally, and by formal asymptotics; see \cite{BM} and references therein.
However, until very recently, there was no rigorous linear (as opposed to
spectral, or normal modes) or nonlinear stability theory
for these waves.

For the physically relevant system \eqref{e:stv},
it turns out that we are able to perform a complete
spectral, linear, and nonlinear
 stability analysis of the associated periodic traveling wave solutions.  Indeed, although
the abstract nonlinear stability theory of \cite{JZ2}, developed for equations of form \eqref{e:conslaw-md} does not
apply directly to the St. Venant equations due to its hyperbolic-parabolic nature and source terms, a suitable
modification of this theory can be made to establish that spectral stability of a given periodic traveling
wave implies nonlinear bounded stability.  This theory will be outlined briefly in the proceeding sections.
The interested reader is invited to find more details in \cite{JZN,N1,N2,NR,BJNRZ}.

The outline of this paper is as follows.  In
the next subsection,
Section \ref{s:ptw},
we briefly review the existence theory for the periodic
traveling wave solutions of the generalized St. Venant equations \eqref{e:stv}.  In particular, we
derive the Hopf bifurcation conditions which guarantee the bifurcation of a family of periodic orbits from
the equilibrium solution.
Then, in Section \ref{s:analytical}, we
outline the known stability theory for the
periodic solutions found in Section \ref{s:ptw}.  The main result of this section is that, under some ``natural" spectral
stability assumptions, the given periodic wave is nonlinearly stable under the PDE dynamics in an appropriate sense.  The
aforementioned spectral stability assumptions are motivated through consideration of the associated Whitham averaged
system and its hyperbolicity, i.e. local well-posedness.  The details of the nonlinear stability proof are
beyond
the scope of the current presentation and hence only an outline of the argument is given; the interested reader is referred to \cite{JZN} for details.

After establishing that (an appropriate sense of) spectral stability implies
nonlinear stability, we then turn our attention in Section \ref{s:spec}
to the verification of the spectral stability assumptions,
restricting our attention to the commonly studied case
$(r,s)=(2,0)$ in \eqref{e:stv}.
We begin by looking in the Hopf and homoclinic limits, as these limits are amenable to direct analysis.
%
%

A striking feature of the
particular example chosen, corresponding to $(r,s)=(2,0)$ in \eqref{e:stv},
is that, in the region of existence of periodic orbits,
{\it all constant solutions are spectrally unstable}
with respect to frequencies
in a neighborhood of the origin.
At a philosophical level, this means that,
for this model, stability can only come about
through dynamical effects having to do with the variation of
the underlying wave, and cannot be understood from a ``frozen-coefficients''
point of view.
To us, this makes it a particularly interesting and illuminating example
from a phenomenological point of view.

At a practical level,
it means that
periodic waves sufficiently close to either the Hopf equilibrium or bounding homoclinic orbit
must be spectrally unstable,
as follows by standard continuity arguments.
Thus, stable periodic waves, if they exist, must be found within
intermediate amplitudes bounded away from either of these
simplifying limits.
%
%

We pursue this line of investigation by conducting a numerical study of the spectrum of the intermediate
amplitude waves found in Section \ref{s:num} and we find
numerically
that there indeed exist periodic waves between the Hopf and
homoclinic orbits which are spectrally stable, and hence are nonlinearly stable by the main theorem in Section \ref{s:analytical}.
We then conclude with a brief conclusion and discussion of the future directions of this project.

\subsection{Periodic Traveling Waves}\label{s:ptw}

We begin with a brief discussion of the traveling wave solutions of the generalized St. Venant equations.
Following \cite{JZN} we restrict ourselves to positive velocities $u>0$ and consider \eqref{e:stv} in Lagrangian
coordinates\footnote{These coordinates are of course equivalent to the Eulerian formulation presented in \eqref{e:stv},
and stability in one formulation is clearly equivalent to stability in the other.  Moreover, we point out
that while the Eulerian formulation is possibly more physically transparent, the Lagrangian formulation is
more analytically convenient.}
\ba \label{e:stvl}
\tau_t - u_x&= 0,\\
u_t+ ((2F)^{-1}\tau^{-2})_x&=
1- \tau^{s+1} u^r +\nu (\tau^{-2}u_x)_x ,
\ea
where $\tau:=h^{-1}$ and now the variable $x$ denotes a Lagrangian marker, rather than a physical location.  Notice
that since the equation \eqref{e:stv} models waves propagating down a ramp, there is no loss in enforcing the restriction
$u>0$\footnote{However, we must always remember to discard any spurious solutions for which $u$ may become negative.}.
In this coordinate frame, a traveling wave solution of \eqref{e:stvl} is a solution
which is stationary in an appropriate moving coordinate frame of the form $x-st$, where
$s\in\RM$ is the wavespeed.  That is, they take the form
\begin{equation}
U(x,t)=\bar U(x-st),
\notag \end{equation}
where $\bar U(\cdot)=(\bar\tau(\cdot),\bar u(\cdot))$ is a solution of the ODE
\be \label{e:secondorderprofile}
\begin{aligned}
-c\bar\tau' - \bar u'&= 0,\\
-c \bar u'+ ((2F)^{-1}\bar\tau^{-2})'&=
1- \bar\tau^{s+1} \bar u^r +\nu (\bar\tau^{-2}\bar u')' .
\end{aligned}
\ee
Integration of the first equation yields $\bar u=\bar u(\tau;q,s):=q-c\bar\tau$, where $q$ is the
corresponding integration constant.  Substitution of this identity into the second
equation yields the second-order scalar profile equation for the function $\bar\tau$:
\be \label{e:profile}
c^2 \bar\tau'+ ((2F)^{-1}\bar\tau^{-2})'=
1- \bar\tau^{s+1} (q-c\bar\tau)^r -c\nu (\bar\tau^{-2}\bar\tau ')' .
\ee
The orbits of \eqref{e:profile} can be studied by simple phase plane analysis.  In particular,
we find that the equilibrium solutions $\tau_0$ satisfy the algebraic
identity
\[
\tau_0^{s+1}(q-c\tau_0)^r=1.
\]

Considering homoclinic solutions then, as in \cite{BJRZ}, under an appropriate normalization we can
assume $\tau_0=1$.  Here, however, we are interested in the periodic orbits of the profile ODE \eqref{e:profile}
which are easily seen not to exist in the case $c=0$.  Indeed, in that case \eqref{e:profile} reduces
to the first order scalar ODE
\be \label{e:zerocprofile}
  \bar\tau'= F\bar\tau^3(\bar\tau^{s+1} q^r-1) ,
\ee
which clearly has no nontrivial solutions with $\bar\tau>0$.  The existence of periodic orbits for non-zero values of $c$
was considered in \cite{N1,N2}, and are generically seen to emerge from a Hopf bifurcation from the equilibrium state
generating a family of periodic orbits which terminate into the bounding homoclinic.  The conditions
for a Hopf bifurcation to occur can be derived from straightforward Fourier analysis as in \cite{BJRZ}.
Indeed, simply notice that the linearization of the profile ODE \eqref{e:profile} about an equilibrium solution
$\tau_0$ (with $q=u_0+c\tau_0)$ is given by
\[
\left(\frac{s+1}{\tau_0}-\frac{cr}{u_0}\right)\tau+(c^2-c_s^2)\tau'+\frac{c\nu\tau''}{\tau_0^2}=0,
\]
where we have used the relation $\tau_0^{s+1}u_0^r=1$.
Taking the Fourier transform, it follows that the Fourier frequency $k$ must satisfy the polynomial equation
\[
\frac{s+1}{\tau_0}-\frac{cr}{u_0}+ik(c^2-c_s^2)-\frac{c\nu k^2}{\tau_0^2}=0.
\]
Evidently, such a $k\in\RM$ exists if and only if $c=c_s$ and
$\left(\frac{s+1}{r}\right)\frac{u_0}{\tau_0}=\left(\frac{s+1}{r}\right)\tau_0^{-(r+s+1)/r}>c_s$, in which case the
solutions are $\pm k_H$ with $k_H\neq 0$.  These translate then to the {\it Hopf bifurcation conditions}
\begin{equation}\label{e:hopf}
c=c_s=\frac{\tau_0^{-3/2}}{\sqrt{F}}\quad\textrm{and}\quad \left(\frac{s+1}{r}\right)\tau_0^{-(r+s+1)/r}>c_s.
\end{equation}
When $(r,s)=(2,0)$, which is the case considered in \cite{JZN}, this
reduces to $F>4$.  This is in agreement with the experiments of \cite{N2,BJNRZ}, which indicate that
when $(r,s)=(2,0)$ and $F>4$ there exists a smooth family of periodic orbits of \eqref{e:profile} parametrized
by the period, which increase in amplitude as the period is increased and finally approaching a limiting homclinic
orbit as the period tends to infinity.

As far as the general existence theory is concerned, we notice that periodic orbits of \eqref{e:profile} correspond
to values $(X,c,q,b)\in\RM^5$, where $X$, $c$, and $q$ denote the period, constant of integration, and wave speed, respectively,
and $b=(b_1,b_2)$ denotes the initial values of $(\tau,\tau')$ at $x=0$ or, equivalently, at $x=X$.  Furthermore, in the
spirit of \cite{Se1,OZ3,OZ4,JZ2,JZ3}, we make the following general assumptions:
\begin{itemize}
\item[(H1)] $\bar \tau>0$, so that all terms in \eqref{e:secondorderprofile}
are $C^{K+1}$, $K\ge 4$.
\item[(H2)] The map $H: \,
\R^5  \rightarrow \R^2$	
taking $(X,c,q,b) \mapsto (\tau,\tau')(X,c,b; X)-b$
is full rank at $(\bar{X},\bar c, \bar b)$,
where $(\tau,\tau')(\cdot;\cdot)$ is the solution operator of \eqref{e:profile}.
\end{itemize}
\noindent
By the Implicit Function Theorem, then, conditions (H1)--(H2) imply that the set of periodic solutions
in the vicinity of $\bar U$ form a smooth $3$-dimensional manifold $\{\bar U^\beta(x-\alpha-c(\beta)t)\}$,
with $\alpha\in \RR$, $\beta\in \RR^{2}$.
%
%
%
%

The goal of our analysis is to understand the modulational stability, i.e. the spectral
and nonlinear time evolutionary stability with respect to small {\it localized} initial perturbations.
We begin by considering the Whitham averaged system corresponding to the dynamical version of \eqref{e:secondorderprofile},
which yields a necessary condition for spectral stability.  These considerations lead to a natural set of spectral
stability assumptions, and under these assumptions we outline the recent nonlinear stability theory
developed in \cite{JZN}.  With this theory in place, we numerically study the spectrum of the linearized
operator for various values of the turbulent parameters $(r,s)$.  In particular, we are able to numerically
find a spectrally stable periodic traveling wave solution of \eqref{e:stvl}.  Recall that, up till now, the existence
of such a stable solution was not at all clear from the known examples, for example
found in \cite{OZ1}.

%
%
%
%
%
%
%
%

\section{Analytical Results}\label{s:analytical}

In this section, we review as briefly as possible the known analytical stability and instability results
concerning the periodic traveling wave solutions of \eqref{e:stvl}.  While this theory is similar
to that developed for the parabolic conservation laws \eqref{e:conslaw-md} developed in \cite{OZ4,JZ2,JZ4},
the extension to the current case involves a number of subtle technical issues associated with lack of parabolicity and nonconservative form.  In particular,
the presence of non-divergence source terms in \eqref{e:stvl} requires a more detailed analysis of the associated Green function
since there are no derivatives to enhance decay.  This is handled by an observation made from the structure of this Whitham averaged equations
that of all the modulations the wave can undergo under low-frequency perturbation, modulations in translation dominate.   This serves
as motivation for a decomposition of the linearized solution operator and allows us to prove the time-asymptotic convergence
of the underlying periodic profile to an appropriate modulation of itself.  The precise statement of this result is the main focus of this section.

%

We begin our study by analyzing the spectral stability of a periodic traveling wave solution of \eqref{e:stvl}.  To this end,
notice that \eqref{e:stvl} can be written in the abstract form
\be\label{ab}
U_t +f(U)_x=(B(U)U_x)_x +g(U)
\ee
and linearizing (\ref{e:stvl}) about $\bar{U}(\cdot)$, we obtain
\be \label{e:lin}
v_t = Lv := (\partial_x B\partial_x   -\partial_x A +C) v,
\ee
where the coefficients
\ba\label{coeffs}
A&:= df(\bu) - (dB(\bu) (\cdot) )\bu_x
=\bp -c & -1\\ -\bar \tau^{_-3}(F^{-1}- 2\nu \bar u_x)&-c\ep ,\\
B&:=B(\bu)= \bp 0&0\\0 & \nu \bar \tau^{-2}\ep,
\quad C:=dg(\bar U)=\bp 0 & 0\\-\bar u^2& -2\bar u\bar \tau\ep
\ea
are periodic functions of $x$.  As the underlying solution $\bar{U}$ depends on $x$ only,
equation \eqref{e:lin} is clearly autonomous in time.  Seeking separated solutions
of the form $v(x,t)=e^{\lambda t}v(x)$, it is clear that the stability of $\bar{U}$ requires
a detailed analysis of the linearized operator $L$.  In particular, we say the underlying periodic
traveling wave is spectrally stable provided the linearized operator $L$ has no spectrum
in the unstable right half plane $\Re(\lambda)>0$.  However, the analysis of the spectrum
of $L$ is made exceedingly difficult by the following two facts: first, as the coefficients of $L$ are $X$-periodic,
Floquet theory implies that the spectrum is purely continuous and hence any spectral instability of the
underlying periodic wave must come from the essential spectrum; secondly, by the translation invariance
of \eqref{e:stvl} it is known that $\bar U'$ is an eigenfunction of $L$ corresponding to $\lambda=0$ and hence
the essential spectrum intersect the imaginary axis in at least one point.  The first issue is dealt with
here by conducting a numerical study of the spectrum as opposed to a analytical spectral stability study.

The second issue on the other hand actually gives us a starting point
for our spectral stability study.
Since it is possible that the spectral curve through the origin might pass through to the unstable
half plane, a natural place to begin our study is to analyze the spectrum of the linearized operator in a neighborhood
of the origin $\lambda=0$ in the spectral plane.  Physically, instability/stability in a neighborhood of the origin corresponds
to the underlying wave being spectrally stable to long-wavelength perturbations, i.e. to slow modulations of the traveling
wave profile.  Thus, we can analyze the long-wavelength stability of a periodic traveling wave by using
a well-developed (formal) physical theory for dealing with such stability problems known as Whitham theory.
In the next section, we summarize recent results concerning the application of Whitham theory to the current situation
and its rigorous verification through the use of Evans function techniques.  This will lead to an analytically necessary
condition for spectral stability and hence to a natural set of stability assumptions similar to those proposed by Schnieder in the context
of reaction-diffusion and related pattern-formation systems \cite{S1,S2,S3}.

\subsection{Whitham averaging and spectral instability}\label{s:whitham}

Very recently a necessary condition for the spectral stability of periodic traveling wave solutions of the generalized St. Venant
equations in Eulerian coordinates \eqref{e:stv} has been derived by a novel relation between the Evans function
and the corresponding linearized Whitham averaged system proposed by Serre \cite{Se1}; see \cite{NR} for complete details.
In particular, the authors show that the linearized dispersion relation obtained from the leading order asymptotics of the Evans function
near the origin can be derived formally through a slow modulation (WKB) approximation yielding the Whitham averaged system.
It follows that the formal homogenization procedures introduced by Whitham \cite{W} and Serre \cite{Se1} yields a necessary
condition for the spectral stability of the underlying periodic traveling wave.  Here, we briefly review this procedure
and its implication for spectral instability of the periodic waves of \eqref{e:stv}, and hence of \eqref{e:stvl}.

As a first step, we let $\eps>0$ be a small perturbation parameter and introduce a set of slow-variables $(x,t)=(\eps X,\eps T)$.  In these
slow variables, we search for a solution of \eqref{e:stv} of the form
\[
(h,u)(X,T)=(h^0,u^0)\left(X,T;\frac{\phi(X,T)}{\eps}\right)+\eps(h^1,u^1)\left(X,T;\frac{\phi(X,T)}{\eps}\right)+\mathcal{O}(\eps^2),
\]
where $\RM\ni y\mapsto (h^j,u^j)(X,T;y)$ are unknown $1$-periodic functions.  It follows then that the local period of oscillation
is $\eps/\partial_X\phi$, where we assume the unknown phase \textit{a priori} satisfies $\partial_X\phi\neq 0$.  Plugging this expansion
into rescaled version of \eqref{e:stv} and collecting like powers of $\eps$ results in a hierarchy of consistency conditions
which must hold.  At the lowest order $\mathcal{O}(\eps^{-1})$, we find that the functions $(h^0,u^0)$
satisfy the corresponding rescaled profile ODE with wavespeed $s$ in the variable $\omega y$, where
\[
s=-\frac{\phi_T}{\phi_X}, \quad\textrm{and}\quad \omega=\phi_T.
\]
Furthermore, notice then that $\omega$ denotes the local frequency and $k=\phi_X$ the local wave number of the modulated wave.
It follows then that $(h^0,u^0)$ can be chosen to agree with a given periodic traveling wave solution of \eqref{e:stv} in the
variable $y$.

Continuing, collecting the $\mathcal{O}(1)$ terms yields the mass conservation law
\[
\partial_y\left(\omega(k,\bar q)h^1-ku^1\right)=-\left(\partial_T h^0+\partial_X u^0\right),
\]
which has a solution if and only if the the right hand side has zero spatial average over a period, i.e. if and
only if
\begin{equation}\label{e:w1}
\partial_T\left(M(X,T)\right)+\partial_X\left(cM(X,T)-q\right)=0,
\end{equation}
where $M(X,T):=\int_0^1h^0(X,T,y)dy$ denotes the corresponding mass functional, $c$ the wave speed, and $q=ch^0-u^0$
the corresponding integration constant.  Together with the consistency condition
\begin{equation}\label{e:w2}
\partial_X k(X,T)+\partial_X(k(X,T)c(X,T))=0
\end{equation}
these equations form a closed first order linear system of partial differential equations known
as the Whitham averaged equations\footnote{Notice that these are not the true Whitham equations for \eqref{e:stv}.  Indeed,
the Whitham equations are the inherently nonlinear equations arising at order $\mathcal{O}(1)$ from substituting the
above expansion in the rescaled version of \eqref{e:stv}.  Upon averaging these nonlinear equations over one spatial period,
however, one arrives at the given linear system; hence its naming as the Whitham {\it averaged} equations.}.

The linear system \eqref{e:w1}-\eqref{e:w2} is seen to be of evolutionary type provided that the non-degeneracy condition
\begin{equation}\label{e:nondegenerate}
\partial_q M(X,T)\neq 0,
\end{equation}
holds. This condition is discussed in both the small-amplitude and small-viscosity regimes in \cite{NR}.  Furthermore, the local well-posedness of this linear system is equivalent with its local hyperbolicity, i.e. the fact that the dispersion relation
$\Delta(\lambda,\nu)$
\[
\Delta(\lambda,\nu):=\det\left(\lambda\frac{\partial(k,M)}{\partial(c,q)}-\nu\frac{\partial(kc,cM-q)}{\partial(c,q)}\right)=0,
\]
with $(\lambda,\nu)\in\CM\times i\RM$ and where all arguments are evaluated at the underlying periodic
wave $(h^0,u^0)$, corresponding to the linearization about the underlying wave $(h^0,u^0)$,
has all real roots.  Notice, in particular, that $\Delta(\lambda,\nu)$ is a homogeneous quadratic polynomial in $\lambda$ and $\nu$.

While hyperbolicity of the Whitham averaged system can heuristically be related to its stability to long-wavelength perturbations,
a rigorous proof of this fact has only been recently given in \cite{NR} through the use of the Evans function, which
we briefly recall here.  Writing the linearization of \eqref{e:stv} about a given $X$-periodic traveling wave solution
as the first order system
\[
Y'=A(\lambda)Y,
\]
the Evans function $D(\lambda,\sigma)$ is defined for $(\lambda,\sigma)\in\CM\times S^1$ via
\[
D(\lambda,e^{\nu})=\det\left(\Psi(X;\lambda)-\sigma\mathbb{I}_3\right),
\]
where $\Psi(\cdot;\lambda)$ denotes the fundamental solution matrix, normalized so that $\Psi(0;\lambda)=\mathbb{I}_3$,
evaluated at the period point $X$.  In particular, notice that $\lambda$ belongs to the spectrum of the linearized operator
if and only if $D(\lambda,\sigma)=0$ for some $\sigma\in S^1$, and hence spectral stability of the underlying wave is equivalent
to the condition that $D(\lambda,\sigma)$ does not vanish for any $\Re(\lambda)>0$ and $\sigma\in S^1$.  Notice, however, that
$D(0,1)=0$ by translation invariance, and hence the spectrum must touch the imaginary axis.  Whether or not the periodic
wave is spectrally stable in a neighborhood of the origin is connected to the hyperbolicity of the Whitham averaged system
through the following theorem.

\begin{thm}[Noble \& Rodrigues \cite{NR}]\label{t:nr}
Let $\bar U$ be a periodic traveling wave solution of \eqref{e:stv} such that the non-degeneracy condition
\eqref{e:nondegenerate} holds.   Then in the limit $(\lambda,\nu)\to (0,0)$ the following asymptotic relation holds:
\[
D(\lambda,e^{\nu})=\Gamma\Delta(\lambda,\nu)+\mathcal{O}((|\lambda|+|\nu|)^3)
\]
for some non-zero constant $\Gamma$.
\end{thm}

That is, the dispersion relation $\Delta(\lambda,\nu)$ agrees to leading order with the Evans function in a neighborhood of the origin.
Recalling that $\Delta$ is a homogeneous in the variables $\lambda$ and $\nu$, introduction of the projective
coordinate $z=\frac{\lambda}{\nu}$ reduces the dispersion relation to the quadratic polynomial
\begin{equation}\label{e:wdisp}
\Delta(z,1)=0,
\end{equation}
whose roots $z_1,z_2$ are distinct so long as the corresponding discriminate is non-zero.  Under this assumption, the implicit
function theorem applies in a neighborhood of $z=z_j$, $\kappa=0$ and hence, in terms of the original spectral variable
$\lambda$ there are two spectral branches
\begin{equation}\label{e:speccurve}
\lambda_j=z_j\nu+\mathcal{O}(\nu^2).
\end{equation}
Thus, if the Whitham system is hyperbolic, corresponding to $z_j\in\RM$, then the two spectral branches emerge from the origin tangent
to the imaginary axis.  This is clearly a necessary condition for spectral stability.  On the other hand, failure of hyperbolicity
of the Whitham system implies that the $z_j$ have non-zero imaginary part and hence the corresponding spectral branches must emerge
from the origin and enter the unstable half plane, immediately yielding spectral instability of the underlying wave.
%

Similar results concerning the spectral verification of the Whitham averaged equations have also been derived
in the viscous conservation law setting \cite{OZ3,Se1}.  Notice however, that while hyperbolicity of the Whitham
system is necessary for the spectral stability of a given periodic traveling wave solution of \eqref{e:stv}, it may
not be sufficient.  Indeed, hyperbolicity of the Whitham system is a first order condition, implying agreement of the spectrum
near the origin along lines.  Thus, the Whitham system will be hyperbolic so long as the spectral curve is tangent to the imaginary
axis at the origin, whether or not the spectral curve then proceeds to the stable or unstable half planes.  Nevertheless
these considerations lead us to a natural set of spectral stability assumptions, which in the next section we show
implies nonlinear stability of the underlying wave.

\subsection{Bloch decomposition and spectral stability conditions}

A particularly useful way to analyze the continuous spectrum of the
linearized operator $L$ is to decompose the problem into a continuous
family of eigenvalue problems through the use of a Bloch decomposition.
To this end, a straightforward application of Floquet theory implies
that the $L^2$ spectrum of the linearized operator $L$ is purely continuous
and corresponds to the union of the $L^\infty$ eigenvalues of the operator $L$
taken with boundary conditions $v(x+X)=e^{i\kappa}v(x)$ for all $x\in\RM$, where $\kappa\in[-\pi,\pi)$
is referred to as the Floquet exponent.  In particular, it follows that $\lambda\in\sigma(L)$
if and only if the spatially periodic spectral problem
\begin{equation}\label{e:lin2}
Lv=\lambda v
\end{equation}
admits a uniformly bounded eigenfunction of the form $v(x)=e^{i\xi x}w(x)$, where $w$ is $X$-periodic.
Substitution of this Ansatz into \eqref{e:lin2} motivates the use of the Fourier-Bloch decomposition
of the spectral problem.

Following \cite{G} then, we define a one-parameter family of linear operators, referred to as the Bloch operators,
via
\[
L_{\xi}:=e^{-i\xi x}Le^{i\xi x},\quad\xi\in[-\pi,\pi)
\]
operating on $L^2_{\rm per}([0,X])$, the space of $X$-periodic square integrable functions.  The spectrum of $L$ is then
seen to be given by the union of the spectra of the Bloch-operators.  Furthermore, since the domain $[0,X]$ is compact
the operators $L_{\xi}$, for each fixed $\xi$, have discrete spectrum in $L^2_{\rm per}([0,X])$ and hence,
by continuity of the spectrum, the spectra of $L$ may be described by the union of countable many continuous surfaces.

Continuing,
taking without loss of generality $X=1$,
we recall that any localized function $v\in L^2(\RM)$
admits an inverse Bloch-Fourier representation
\[
v(x)=
\left(\frac{1}{2\pi }\right) \int_{-\pi}^{\pi} e^{i\xi x}\hat v(\xi, x) d\xi
\]
where the functions $\hat v(\xi,\cdot)=\sum_{j\in\ZM}e^{2\pi ijx}\hat v(\xi+2\pi j)$ belongs to $L^2_{\rm per}([0,X])$ for
each $\xi$, where here $\hat v(\cdot)$ denotes with a slight abuse of notation the usual Fourier transform of the function
$v$ in the spatial variable $x$.  By Parseval's identity it is seen that the Bloch-Fourier transformation
$v(x)\to\hat v(\xi,x)$ is an isometry of $L^2(\RM)$, i.e.
\[
\|v\|_{L^2(\RM)}=\int_{-\pi}^\pi\int_0^X|\hat v(\xi,x)|^2dx~d\xi=:\|\hat g\|_{L^2(\xi;L^2(x))}.
\]
Furthermore, this transformation is readily seen to diagonalize the periodic-coefficient linearized operator $L$,
yielding the inverse Bloch-Fourier transform representation
\[
e^{Lt}v(x)=\frac{1}{2\pi}\int_{-\pi}^\pi e^{i\xi x}e^{L_{\xi}t}\hat g(\xi,x)d\xi
\]
effectively relating the behavior of the linearized system to that of the diagonal operator $L_\xi$.

Together with the long-wavelength stability analysis in the previous section, we now state our main
spectral stability assumptions in terms of the diagonal Bloch operators $L_{\xi}$.
\begin{itemize}
\item[(D1)] $\sigma(L_\xi)\subset\{\lambda\in\CM:\Re\lambda<0\}$ for all $\xi\neq 0$.
\item[(D2)] There exists a constant $\theta>0$ such that $\Re\sigma(L_\xi)\leq-\theta|\xi|^2$ for all $|\xi|\ll 1$.
\item[(D3')] $\lambda=0$ is an eigenvalue of $L_0$ of multiplicity exactly two.\footnote{Note that the zero eigenspace of $L_0$,
corresponding to variations along the three-dimensional manifold of periodic solutions in directions for which the
period does not change \cite{Se1,JZ2}, is at least two-dimensional by the linearized existence theory and assumption (H2).}
\end{itemize}
Notice that assumption $(D1)$ implies weak hyperbolicity of the Whitham averaged system, while $(D2)$ corresponds to ``diffusivity"
of the large-time ($\sim$ small frequency) behavior of the linearized operator $L$.  Moreover, $(D3')$ holds generically and can be
directly verified through the use of the Evans function arguments as in \cite{N1}.  Finally, we point out that $(D1)-(D3')$ are balance
law analogues of the spectral assumptions introduced by Schneider for reaction-diffusion equations \cite{S1,S2,S3}.

Furthermore, we make the following non-degeneracy hypothesis:
\begin{itemize}
\item[(H3)] The roots $z_j$ of \eqref{e:wdisp} are distinct.
\item[(H4)] The eigenvalue $0$ of $L_0$ is non-semisimple, i.e. $\dim\ker(L_0)=1$.
\end{itemize}
Conditions $(H1)-(H4)$ generically imply that (D2) hold\footnote{This amounts to nonvanishing $b_j$ in the
Taylor series expansion $\lambda_j(\xi)=-iz_j\xi-b_j\xi^2+o(|\xi|^2)$ guaranteed by Lemma \ref{blochfacts} given $(H1)-(H4)$,
$(D1)$, and $(D3')$.}.  Moreover, (H3) corresponds to strict hyperbolicity of the Whitham averaged system, and implies
the analyticity of the spectrum in a neighborhood of the origin.  Specifically, since $\Delta(0,1)\neq 0$,
as is readily seen in \cite{NR}, it follows that the roots $z_j$
of \eqref{e:wdisp} are non-zero and distinct by (H3) and hence relation \eqref{e:speccurve} and standard spectral
perturbation theory \cite{K} implies the spectral curves $\lambda_j=\lambda_j(\nu)=\lambda_j(i\xi)$
are analytic functions of $\xi$ in a neighborhood of $\xi=0$.  Finally, notice that assumptions $(D3')$ and $(H4)$ imply the existence of a Jordan block at $(\lambda,\xi)=(0,0)$.
In particular, we have the following spectral preparation result.

\begin{lemma}[\cite{JZN}]\label{blochfacts}
Assuming (H1)--(H4), (D1), and (D3'), the eigenvalues
$\lambda_j(\xi)$ of $L_\xi$
are analytic functions and the Jordan structure of the zero
eigenspace of $L_0$ consists
of a $1$-dimensional kernel and a single Jordan chain of height $2$,
where the left kernel of $L_0$ is spanned by the constant
function $\tilde f\equiv (1,0)^T$, and $\bar u'$ spans the right
eigendirection lying at the base of the Jordan chain.
Moreover, for $|\xi|$ sufficiently small,
there exist right and left eigenfunctions
$q_j(\xi, \cdot)$ and $\tilde q_j(\xi, \cdot)$
of $L_\xi$ associated with $\lambda_j(\xi)$ which
are analytic in $\xi$ for in a neighborhood of $\xi=0$.
Furthermore, $\langle \tilde q_j,q_k\rangle= \delta_j^k$%
%
\end{lemma}

\br
Notice that the results of Lemma \ref{blochfacts} are somewhat unexpected since, in general, eigenvalues
bifurcating from a non-trivial Jordan block typically do so in a nonanalytic fashion,
rather being expressed in a Puiseux series in fractional powers of $\xi$.  The
fact that analyticity prevails in our situation is a consequence of the very special structure
of the left and right generalized null-spaces of the unperturbed operator $L_0$, and the special
forms of the equations considered.
\er

From the standpoint of obtaining a nonlinear stability result, the existence of the non-trivial Jordan block
over the translation mode suggests that one can not expect
traditional orbital asymptotic stability of the original periodic traveling wave in any standard
$L^p$ or $H^s$ norm; see \cite{OZ2}.  Nevertheless, following the
ideas of \cite{JZ2} we are able to prove nonlinear asymptotic stability to an appropriate
modulation of the original wave, and hence an $L^\infty$ stability result for the underlying wave.
The technical details driving this nonlinear stability argument are
beyond the scope
of what we wish to discuss here; the interested reader can see \cite{JZ2,JZN}.
However, for a sense of completeness we recall here the general outline of the argument in the next section.

\subsection{Nonlinear stability: a guided tour}

Here, we wish to recall the basic ideas behind the nonlinear stability of periodic traveling
wave solutions of the St. Venant equations \eqref{e:stv}.  For technical reasons, we find it essential
to utilize the Lagrangian formulation \eqref{e:stvl} throughout this analysis.  To begin, let $\bar{U}(x)$
denote a periodic traveling wave solution of \eqref{e:secondorderprofile} and let $\tilde U(x,t)$ denote
any other solution of \eqref{e:stvl}.  Our goal is to prove that if $\tilde U(x,0)$ is sufficiently close
to $\bar U(x)$ in a suitable norm, then it remains close for all future times $t>0$.  To this end,
define the nonlinear perturbation variable
\begin{equation}\label{e:res}
v(x):=\tilde U(x+\psi(x,t))-\bar U(x),
\end{equation}
where 
$\psi:\RM^2\to\RM$ is a modulation function to be chosen later.
Our starting point is the following observation: by a direct computation and Taylor expansion,
the nonlinear residual \eqref{e:res} is seen to satisfy
\[
\left(\partial_t-L\right)v=\left(\partial_t-L\right)\bar U'(x)\psi-Q_x+T+P+R_x+\partial_t S,
\]
where
\begin{align*}
P&=(0,1)^T O(|v|(|\psi_{xt}|+|\psi_{xx}|+|\psi_{xxx}|),\\
Q&:=f(\tilde{U}(x+\psi(x,t),t))-f(\bar{U}(x))-df(\bar{U}(x))v=\mathcal{O}(|v|^2),\\
T&:=
(0,1)^T\Big((\tilde{U}(x+\psi(x,t),t))^2-(\bar{U}(x))^2-\bar{U}(x))v\Big)=
(0,1)^T\mathcal{O}(|v|^2),\\
R&:= v\psi_t + v\psi_{xx}+  (\bar U_x +v_x)\frac{\psi_x^2}{1+\psi_x},\quad\textrm{and}\\
S
&=O(|v|(|\psi_x|).
\end{align*}
Letting $G(x,t;y)$ denote the Green function of \eqref{e:lin2} then, an application of Duhamel's formula
implies the nonlinear residual must satisfy the integral equation
\ba\label{prelim}
  v(x,t)&=\psi (x,t) \bar U'(x)+\int^\infty_{-\infty}G(x,t;y)v_0(y)\,dy  \\
  &\quad
  + \int^t_0 \int^\infty_{-\infty} G(x,t-s;y)
  (-Q_y+T+ R_y + S_t ) (y,s)\,dy\,ds.
\ea
In order to obtain control over $v$ in a given $H^s$ norm then, we seek to obtain pointwise
bounds on the Green function $G$ and an appropriate expression for $\psi$ for which
\eqref{prelim} becomes susceptible to an iteration argument.  Since, as expected, the low-frequency
behavior of the solution operator near the neutral eigenvalue $\lambda=0$ is the most difficult
to control, we decompose the solution operator $S(t)=e^{Lt}$ corresponding to the linearized operator $L$
into high and low-frequency components.

To this end, standard spectral perturbation theory \cite{K} implies that the total eigenprojection $P(\xi)$
onto the eigenspace of $L_{\xi}$ associated with the eigenvalues $\lambda_j(\xi)$ described in \eqref{e:speccurve}
is well-defined and analytic in $\xi$ for $|\xi|$ sufficiently small, since these (by discreteness of the
spectra of $L_{\xi}$) are separated at $\xi=0$ from the rest of the spectrum of $L_0$.
Choosing an appropriate cut-off function $\phi(\xi)$ supported
in a small neighborhood of the origin and identically one in a slightly smaller neighborhood,
we split $S(t)$ into a low-frequency part
\[
S^I(t)u_0:=\frac{1}{2\pi}\int_{-\pi}^\pi e^{i\xi x}\phi(\xi)P(\xi)e^{L_{\xi}t}\hat u_0(\xi,x)d\xi
\]
and the associated high-frequency part
\[
S^{II}(t)u_0:=\frac{1}{2\pi}\int_{-\pi}^\pi e^{i\xi x}(1-\phi(\xi))P(\xi)e^{L_{\xi}t}\hat u_0(\xi,x)d\xi.
\]

We begin by analyzing $S^{II}$.  Fairly routine semigroup estimates \cite{He,Pa} imply the bounds
\begin{align*}
\left\|e^{L_\xi t}g\right\|_{L^2([0,X])}&\lesssim e^{-\theta t}\|g\|_{L^2([0,X])}\\
\left\|\partial_x e^{L_\xi t}g\right\|_{L^2([0,X])}&\lesssim e^{-\theta t}t^{-1/2}\|g\|_{L^2([0,X])}\\
\left\| e^{L_\xi t}\partial_xg\right\|_{L^2([0,X])}&\lesssim e^{-\theta t}t^{-1/2}\|g\|_{L^2([0,X])}
\end{align*}
for all times $t>0$ and some constant $\theta>0$.  Since the Bloch-Fourier transform is an isometry of $L^2$, it follows
by standard $L^p$ interpolation and Sobolev's inequality that there exists a constant $\theta>0$ such that
\[
\|S^{II}(t)\partial_x^l g\|_{L^p(\RM)}\lesssim e^{-\theta t}t^{-l/t}\|g\|_{L^2(\RM)}.
\]
for all $0\leq l\leq 2$ and, similarly, one obtains an estimate on terms of the form $\|S^{II}(t)\partial_t^m g\|_{L^2(\RM)}$.

Next, we seek analogous bounds on the low-frequency portion of the
solution operator.  This however is complicated by the presence
of spectral curves which touch the imaginary axis at the origin, and hence a more delicate analysis is necessary.  We begin by
denoting by
\[
G^I(x,t;y):=S^I(t)\delta_y(x)
\]
the Green kernel associated with $S^{I}$.  Furthermore, recalling Lemma \ref{blochfacts}, for $|\xi|$ sufficiently small
we denote by $q_j(x,\xi)$ and $\tilde q_j(x,\xi)$ the right and left eigenfunctions
of the Bloch operator $L_\xi$, respectively, associated with the spectral curves $\lambda_j(\xi)$ bifurcating from the $(\xi,\lambda_j(\xi))=(0,0)$
state and we enforce the normalization condition $\left<\tilde q_j(\cdot,\xi),q_k(\cdot,\xi)\right>_{L^2([0,X])}=\delta_j^k$.
It follows then that we can express the low-frequency Green function as
\[
G^I(x,t;y)=\left(\frac{1}{2\pi}\right)\int_\RM e^{i\xi(x-y)}\phi(\xi)\sum_{j=1}^2e^{\lambda_j(\xi)t}q_j(\xi,x)\tilde q_j(\xi,y)^* d\xi,
\]
where $*$ denotes the complex conjugate transpose.  Notice that this Bloch expansion for the Green kernel is analogous to using a Fourier
representation in the constant coefficient case.
%
Similary as in the constant--coefficient case, we may read off decay
from the spectral representation using the following generalization
of the Hausdorff--Young inequality.

\bl[Generalized Hausdorff--Young inequality \cite{JZ2}]\label{l:hy}
\be\label{e:hy}
\|u\|_{L^p(x)}\leq\|\hat u\|_{L^q(\xi;L^p(0,X))},\quad\textrm{ for }q\leq 2\leq p\quad\textrm{ and }\frac{1}{p}+\frac{1}{q}=1
\ee
\el

\begin{proof}
Relation \eqref{e:hy}
holds in the extremal cases $p=2$ and $p=\infty$ by Parseval's identity and the triangle inequality, respectively.  This generalized
version of the Hausdorff--Young inequality then holds for all stated pairs $(p,q)$ by a generalized version of Riesz-Thorin interpolation theorem;
see Appendix A of \cite{JZ2} for more details.
\end{proof}

In order to analyze the decay properties of the kernel $G^I$ in $t$, we notice
that while assumption (D2) implies $e^{\lambda_j(\xi) t}\lesssim e^{-\theta |\xi|^2t}$ for all $t>0$, this does not immediately
yield decay since the presence of the Jordan block (guaranteed by assumptions (D3') and (H4)) implies
$q_1(\xi)\tilde q_1(\xi)\sim\xi^{-1}$, where $q_1(0)(x)=\bar U'(x)$ corresponds to the translation mode.
%
Following this intuition, we find that
the terms in the Bloch expansion of the Green kernel $G^I$ not
associated with the Jordan block near $\xi=0$ decay in $L^p(x)$ as $\|e^{-\theta\xi^2 t}\|_{L^q(\xi)}\lesssim t^{-\frac{1}{2}(1-1/p)}$,
i.e. at the rate of a heat kernel, while the portion of the kernel associated with the Jordan block decays in $L^\infty(x)$
as $\|\xi^{-1}e^{-\theta\xi^2t}\|_{L^1(\xi)}\lesssim 1$.  These bounds, derived from the Bloch-norm Hausdorff--Young inequality
discussed above should be compared with
those found by
weighted-energy estimate methods of Schneider \cite{S1}.

Using the above $L_x^p\to L_\xi^qL_z^p$ bounds then it follows that the low-frequency Green kernel $G^I$ can be decomposed
as
\[
G^I(x,y;t)=\bar U'(x)e(x,t;y)+\tilde G^I(x,t;y)
\]
where the residual $\tilde G^I$ and amplitude $e$ satisfy the bounds
\begin{align*}
\sup_{y\in\RM}\|\tilde G^I(\cdot,t;y)\|_{L^p(x)}&\lesssim(1+t)^{-\frac{1}{2}(1-1/p)}\\
\sup_{y\in\RM}\|e(\cdot,t;y)\|_{L^p(x)}&\lesssim(1+t)^{-\frac{1}{2}(1-1/p)}
\end{align*}
for all $2\leq p\leq\infty$.  Furthermore, it can be shown that the derivatives of these functions decay in time even faster
according to the variable and order of differentiation.
These estimate can be used to control the ``free-evolution" type terms
appearing in the integral equation \eqref{prelim}.  To control the integrals associated with the (implicit) source terms,
arguments like those outlined above can be used to obtain the estimates
\begin{align*}
\left\|\int_{\RM}\tilde \partial_y^rG(\cdot,t;y)f(y)dy\right\|_{L^p(x)}&\lesssim(1+t)^{-\frac{1}{2}(1/q-1/p)-r/2}\|f\|_{L^q\cap H^1}\\
\left\|\int_\RM \partial_{x,y,t}e(\cdot,t;y)f(y)dy\right\|_{L^p(x)}&\lesssim(1+t)^{-\frac{1}{2}(1/q-1/p)}\|f\|_{L^q}
\end{align*}
where $1\leq q\leq 2\leq p\leq\infty$ and $r=0,1$.

With these preparations, we return to \eqref{prelim} and define
\[
\psi(x,t):=-\int_0^t\int_{\RM}e(x,t-s;y)\left(-Q_y+R_y+\left(\partial_t+\partial_y^2\right)S\right)dy~ds
\]
and note that this choice cancels the ``bad" term $\bar U'(x)e$ in the decomposition $G^I=\bar U'e+\tilde G^I$.
Furthermore, using \eqref{prelim} this choice results in a closed system in the variables $(v,\psi_x,\psi_t)$,
where now $v$ satisfies
\[
v(x,t)=\int_0^t\int_{\RM}\tilde G^I(x,t-s;y)\left(-Q_y+R_y+\left(\partial_t+\partial_y^2\right)S\right)dy~ds,
\]
and
\[
\psi_{x,t}(x,t)=\int_0^t\int_{\RM}e_{x,t}(x,t-s;y)\left(-Q_y+R_y+\left(\partial_t+\partial_y^2\right)S\right)dy~ds.
\]
Recalling then that $(Q,R,S)=\mathcal{O}(|v,\psi_{x,t}|^2)$ and $\tilde G^I$ and $e_{x,t}$ decay at Gaussian rates in $L^p(x)$,
nonlinear stability follows by a direct iteration (contraction mapping) argument as in the more well
familiar viscous shock case (where the profile decays to constant solutions).  As with the linearized bounds
derived above, these cancelation computations, carried out in the physical variables, should be compared
to the cancelation computations carried out in the frequency domain by Schneider for the reaction diffusion case.
In particular, we arrive at the main theorem of \cite{JZN}.

\begin{thm}\label{t:nl}
Assuming (H1)--(H4) and (D1)--(\DDD),
let $\bar U=(\bar \tau, \bar u)$ be a traveling-wave solution
of \eqref{e:stvl} satisfying
the derivative condition
\be\label{froudebd}
\nu \bar u_x < F^{-1}.
\ee
Then, for some $C>0$ and $\psi \in W^{K,\infty}(x,t)$, $K$ as in (H1),
\ba\label{eq:smallsest}
\|\tilde U-\bar U(\cdot -\psi-ct)\|_{L^p}(t)&\le
C(1+t)^{-\frac{1}{2}(1-1/p)}
\|\tilde U-\bar U\|_{L^1\cap H^K}|_{t=0},\\
\|\tilde U-\bar U(\cdot -\psi-ct)\|_{H^K}(t)&\le
C(1+t)^{-\frac{1}{4}}
\|\tilde U-\bar U\|_{L^1\cap H^K}|_{t=0},\\
\|(\psi_t,\psi_x)\|_{W^{K+1,p}}&\le
C(1+t)^{-\frac{1}{2}(1-1/p)}
\|\tilde U-\bar U\|_{L^1\cap H^K}|_{t=0},\\
\ea
and
\ba\label{eq:stab}
\|\tilde U-\bar U(\cdot-ct)\|_{ L^\infty}(t), \; \|\psi(t)\|_{L^\infty}&\le
C
\|\tilde U-\bar U\|_{L^1\cap H^K}|_{t=0}
\ea
for all $t\ge 0$, $p\ge 2$,
for solutions $\tilde U$ of \eqref{e:stvl} with
$\|\tilde U-\bar U\|_{L^1\cap H^K}|_{t=0}$ sufficiently small.
In particular, $\bar U$ is nonlinearly bounded
$L^1\cap H^K\to L^\infty$ stable.
\end{thm}

Theorem \ref{t:nl} asserts not only bounded $L^1\cap H^K\to L^\infty$ stability, a very weak notion of stability,
but also asymptotic convergence of $\tilde U$ to the space-time modulated wave $\bar U(x-\psi(x,t))$.
The fact that we don't get decay to $\bar U$ is the fundamental reason why we had to factor out
the translation mode from the low-frequency Green kernel $G^I$: indeed, the nonlinear (source)
terms in \eqref{prelim} can not be considered as asymptotically negligible without asymptotic decay of the
modulation $\psi$ to zero.  After factoring out the translation mode, however, it is found that the perturbation
$v$ decays to zero as $\psi_x$, which, thanks to the diffusive nature of the linearized operator, decays
at a rate $t^{1/2}$ (in $L^2$) faster than $\psi$.  Finally,
we also note that the derivative condition \eqref{froudebd} is effectively an upper bound on the amplitude of the underlying
periodic wave and is a technical condition needed to control the $H^K$ norm of the perturbation $v$ in terms of the $L^2$ norm;
this technical detail is beyond the scope of our current presentation, and interested readers are referred to \cite{JZN}.
It should be noted, however, that \eqref{froudebd} is satisfied when either the wave amplitude or viscosity coefficient $\nu$
is sufficiently small, and is seen to be satisfied for \emph{all} roll-waves computed numerically in \cite{N2} and here
in Section \ref{s:num}.


Our next goal is to verify, at least numerically, the spectral stability conditions necessary in Theorem \ref{t:nl}
to conclude nonlinear (bounded) stability of the underlying periodic wave $\bar U$.

\section{Spectral Stability}\label{s:spec}

\begin{figure}[htbp]
\begin{center}
$\begin{array}{cc}
(a) \includegraphics[trim=0 0 0 0, scale=.35]{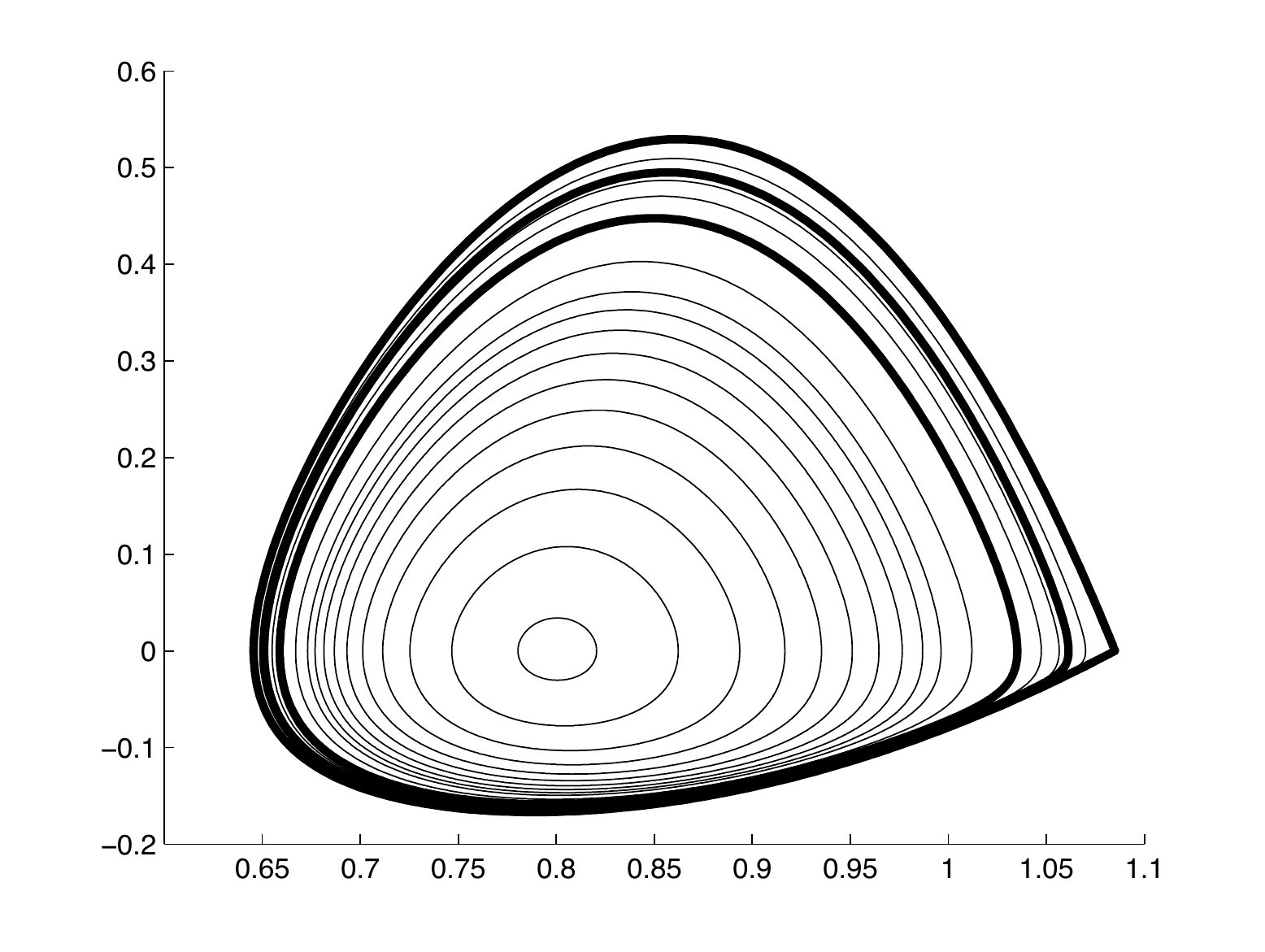} & (b) \includegraphics[trim=0 0 0 0, scale=.35]{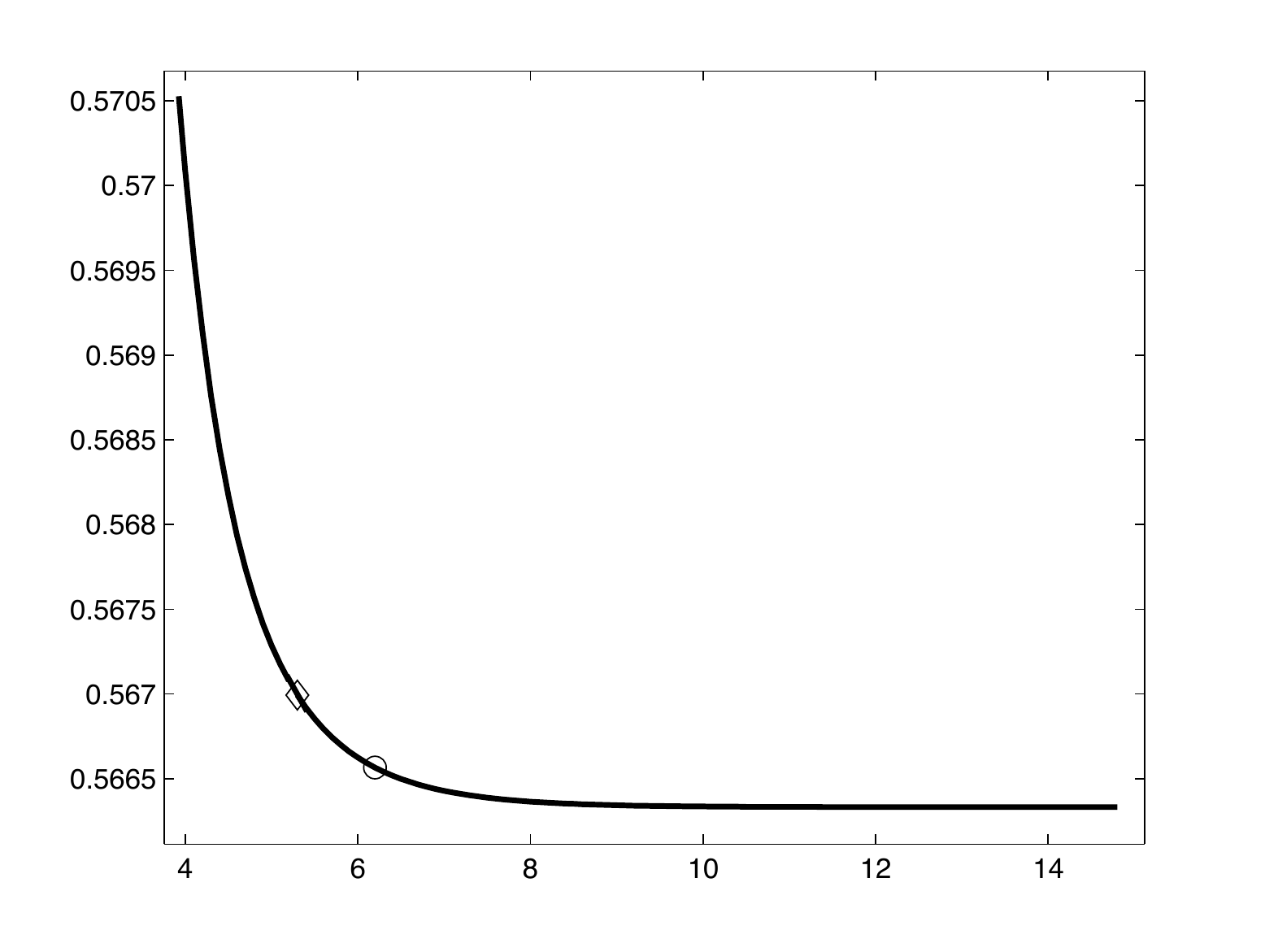}
\end{array}$
\end{center}
\caption{(a) A typical phase portrait depicting a family of periodic orbits parameterized by the wave speed $c$ generated through a Hopf bifurcation
at the enclosed equilibrium solution.  The inner and outer most bold orbits correspond to the lower (period $X\approx 5.3$) and upper ($X\approx 20.6$) stability
boundaries, while the bold orbit in between corresponds to the stable periodic traveling wave solution ($X\approx 6.2$) depicted in Figure \ref{f:stable}.
  (b) A plot of the period $X$ versus the wavespeed $c$ of the corresponding periodic traveling wave.
Notice that all periodic orbits sufficiently close to the bounding
homoclinic have approximately the same wavespeed and hence locally resemble in both shape and speed the limiting homoclinic wave.
The diamond
denotes the lower stability boundary, while the circle signifies the stable solution depicted in Figure \ref{f:stable}.
}
\label{f:phase}
\end{figure}

Now that we have established that the spectral stability of a given periodic traveling wave solution of \eqref{e:stvl}
implies nonlinear stability in the sense of Theorem \ref{t:nl}, we continue our investigation by analyzing
the spectral stability question.  To begin, we restrict to the commonly studied case $(r,s)=(2,0)$
and depict in Figure \ref{f:phase}(a) a typical phase portrait for the corresponding profile ODE \eqref{e:profile}
in the $\tau,\tau'$ variables.  In this case, all periodic orbits arise through a Hopf bifurcation, corresponding
to minimun period $X\approx 3.9$, as $c$ is decreased through
the critical wavespeed $c_s$, as depicted in Figure \ref{f:phase}(b).
An interesting feature here is that the upper stability boundary, corresponding to the
bold orbit of greatest amplitude, is nearly indistinguishable from the limiting homoclinic orbit in both shape, by Figure \ref{f:phase}(a),
and speed, by Figure \ref{f:phase}(b).

The goal of this section is to attempt to find a spectrally stable solution of the St. Venant equation
\begin{equation}
\begin{aligned}\label{e:stvcommon}
\tau_t - u_x&= 0,\\
u_t+ ((2F)^{-1}\tau^{-2})_x&=
1- \tau u^2 +\nu (\tau^{-2}u_x)_x ,
\end{aligned}
\end{equation}
considered here again in Lagrangian coordinates,
There are two natural limits in which the spectrum of the corresponding linearized operator $L$ seems amenable to direct analysis.
The first is a small-amplitude, i.e. Hopf, limit as one approaches the enclosed equilibrium solution, while the other corresponds
to the large-period limit as the periodic wave approaches the bounding homoclinic in phase space.  In the next section,
we recall recent results of \cite{BJRZ} concerning stability in these distinguished limits.

\subsection{Hopf and homoclinic limits}

\begin{figure}[htbp]
\begin{center}
$\begin{array}{cc}
(a) \includegraphics[trim=0 0 0 0, scale=.25]{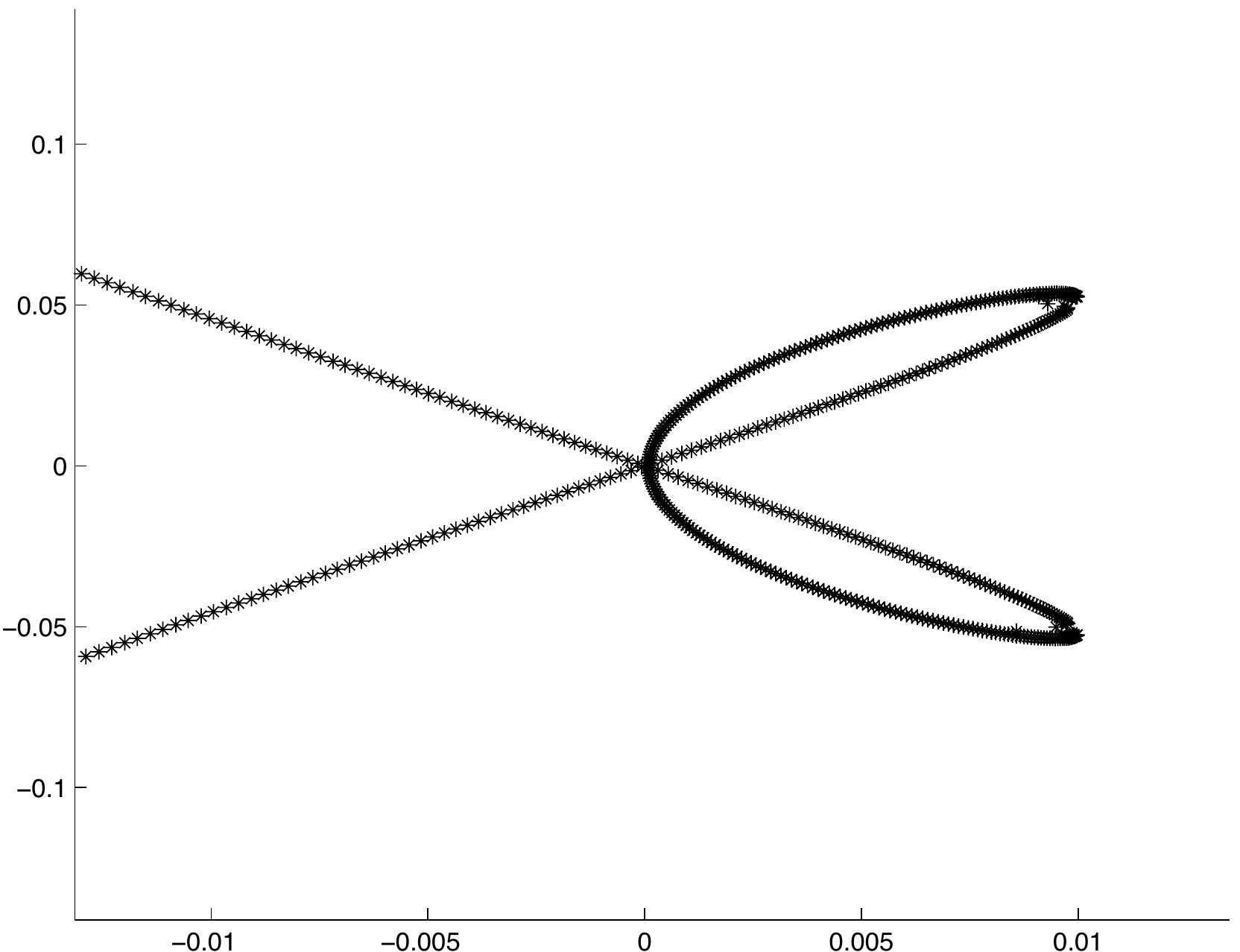}& (b) \includegraphics[trim=0 0 0 0, scale=.2]{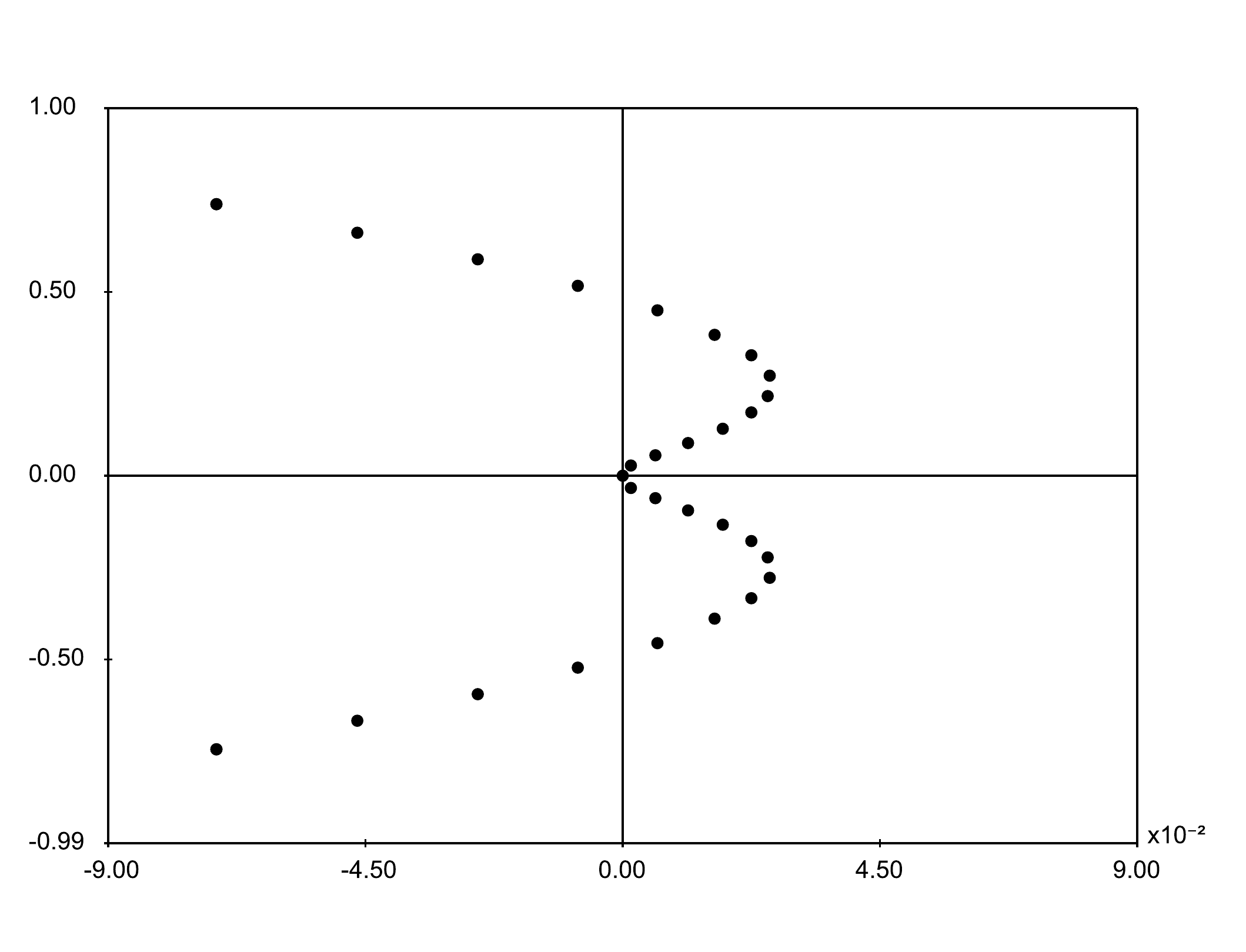}
\end{array}$
\end{center}
\caption{(a) The essential spectrum for the unstable constant solution at the Hopf bifurcation point is depicted.
 (b) The essential spectrum of the bounding homoclinic is shown.
 %
Both of these spectral plots, as well as the rest presented throughout this paper,
are of $\Re(\lambda)$ vs. $\Im(\lambda)$ and were generated using the SpectrUW package developed at the University of Washington \cite{DK},
 which is designed to find the essential spectrum of linear operators with periodic coefficients by using Fourier-Bloch decompositions
 and Galerkin truncation; see \cite{CuD,CDKK,DK} for further information and details concerning convergence.}
 \label{f:cst}
\end{figure}

We begin our search for stable periodic traveling wave solutions of \eqref{e:stvcommon} by considering the
small amplitude limit in which one approaches the enclosed equilibrium solution.  More generally, we consider
the stability of the equilibrium solutions, which satisfy the relation $\tau_0u_0^2=1$.  
To this end, we recall from \cite{BJRZ} that the linearization of \eqref{e:stvcommon} about an equilibrium solution $(\tau_0,u_0)$
satisfying $\tau_0u_0^2=1$ has as an associated dispersion relation between the eigenvalue $\lambda$ of the linearized
operator $L$ and the Fourier frequency $k$ given by
\[
\lambda^2 +
\left[\frac{r}{u_0}-2ick+\frac{\nu k^2}{\tau_0^2}\right]\lambda+
ik\left[\frac{s+1}{\tau_0}-\frac{cr}{u_0}+ik (c^2-c_s^2)-\frac{c\nu k^2}{\tau_0^2}\right]=0.
\]
Notice that the eigenvalues corresponding to frequency $k=0$ are $\lambda=-u_0^{-1}$, which remains negative for $|k|\ll 1$,
and $\lambda=0$.  To determine the behavior of the $\lambda=0$ for small nonzero $k$, we
Taylor expand the dispersion relation about $(\lambda,k)=(0,0)$ with $\lambda=\lambda(k)$ and find that the
spectral curve $\lambda(k)$ must satisfy
\[
\lambda'(0)=-i\left[\frac{u_0}{2\tau_0}-c\right],
\]
indicating stability, while
\[
\frac{1}{2}\lambda''(0)=\frac{u_0}{2}\left[\left(i\lambda'(0)+c\right)^2-c_s^2\right]=\frac{u_0}{2}\left[\left(\frac{u_0}{2\tau_0}\right)^2-c_s^2\right].
\]
Recalling the Hopf bifurcation conditions \eqref{e:hopf},
we find
that if the equilibrium solution $\tau_0$ corresponds
to a Hopf bifurcation point of the profile ODE \eqref{e:profile}
then we must have $\lambda''(0)>0$, yielding
instability of the equilibrium (Hopf) point.  In particular, we see that in the regime of existence of periodic waves, i.e. $F>4$,
\emph{all} constant solutions are spectrally unstable.  Therefore, since the spectrum of the linearized operator $L$
changes continuously as we nearby periodic orbits we conclude that all periodic
traveling wave solutions of \eqref{e:stvcommon} solutions must be spectrally unstable in the small amplitude limit.
This is verified numerically in Figure \ref{f:cst}(a).

Next, we turn to the large-period limit as the periodic orbits approach the bounding homoclinic profile in phase space.
In this case, we can use the same arguments as in the Hopf limit in order to determine the stability of the limiting
end state of the homoclinic orbit, which recall determines the essential spectrum of the homoclinic \cite{He}.
It follows that the limiting endstate is spectrally unstable in a neighborhood of the origin whenever the Hopf bifurcation
condition $F>4$ hold, and this is numerically verified in Figure \ref{f:cst}(b).  Therefore, we should
not expect to find any spectrally stable periodic waves in the homoclinic limit.

\begin{figure}[htbp]
\begin{center}
$
\begin{array}{ccc}
  \includegraphics[scale=.25]{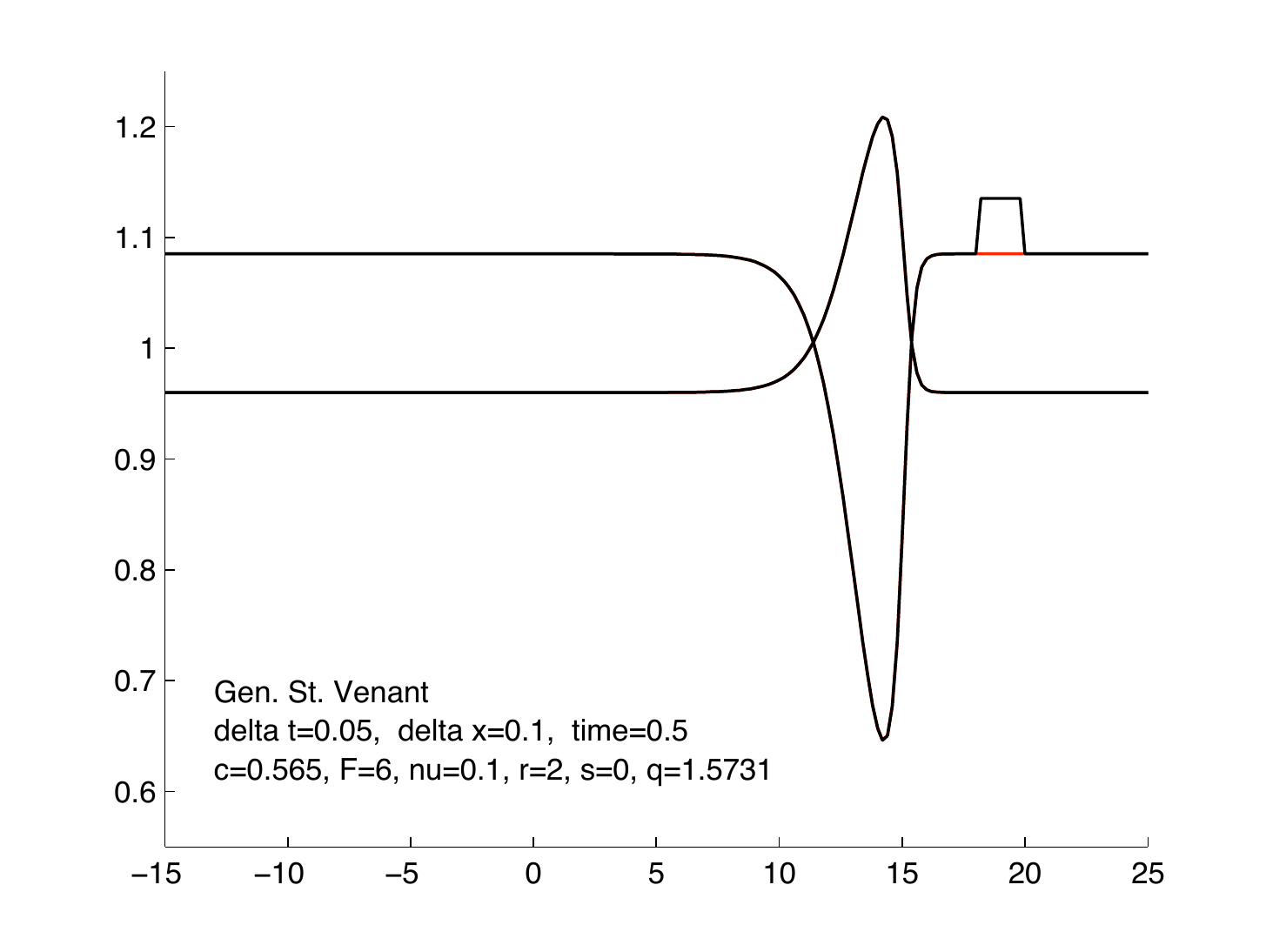}&\includegraphics[scale=0.25]{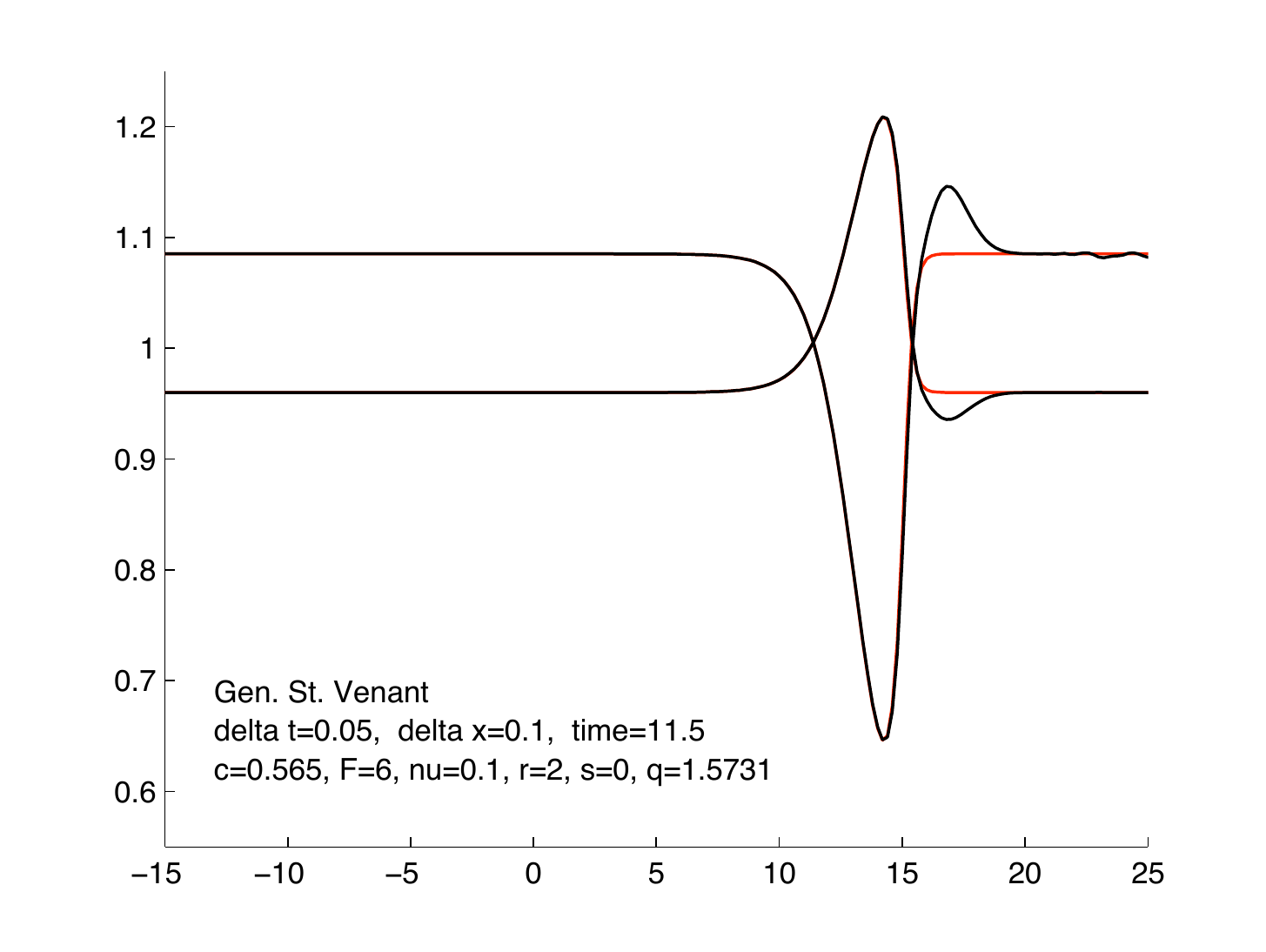}& \includegraphics[scale=.25]{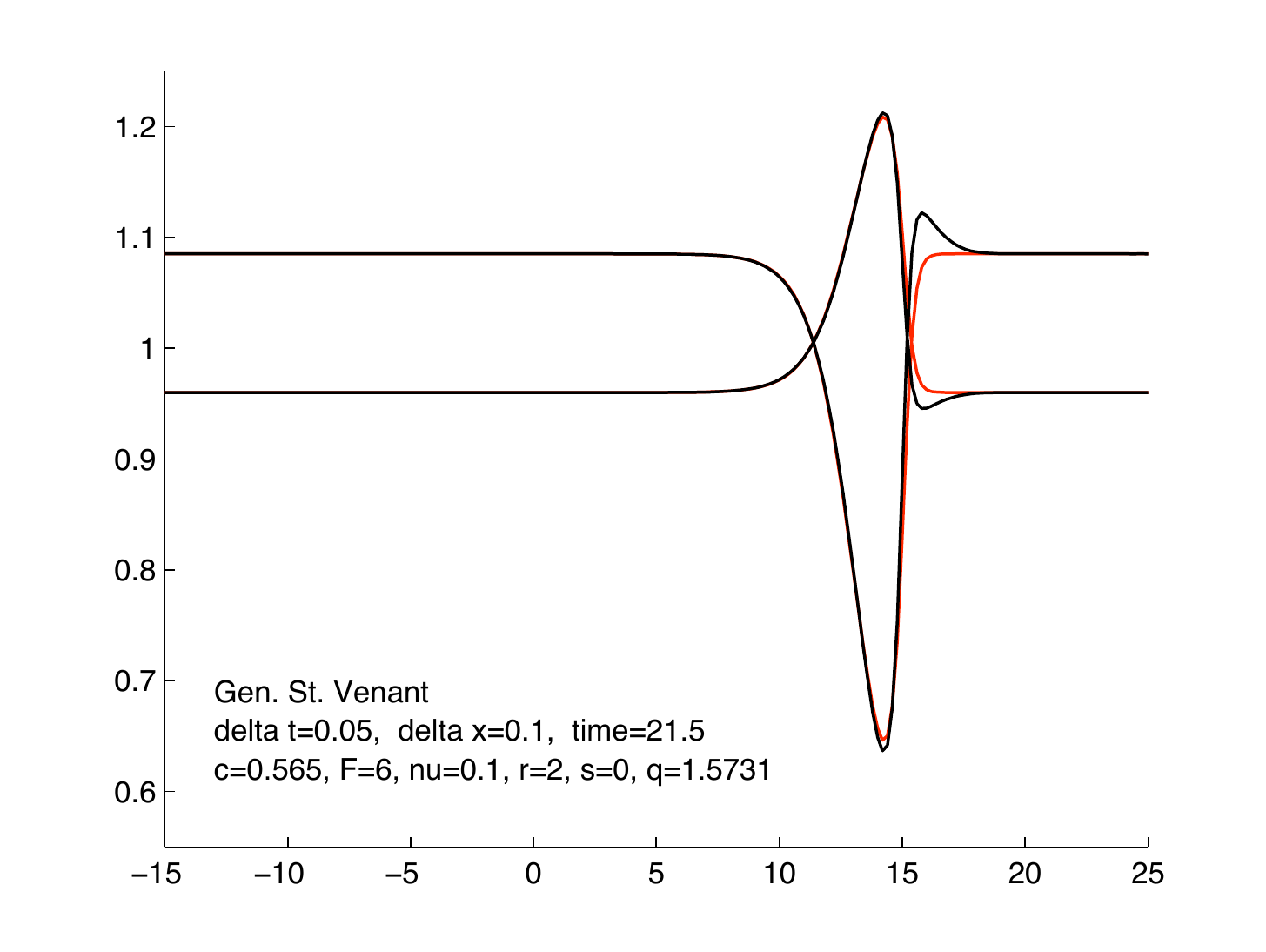}\\
  \includegraphics[scale=.25]{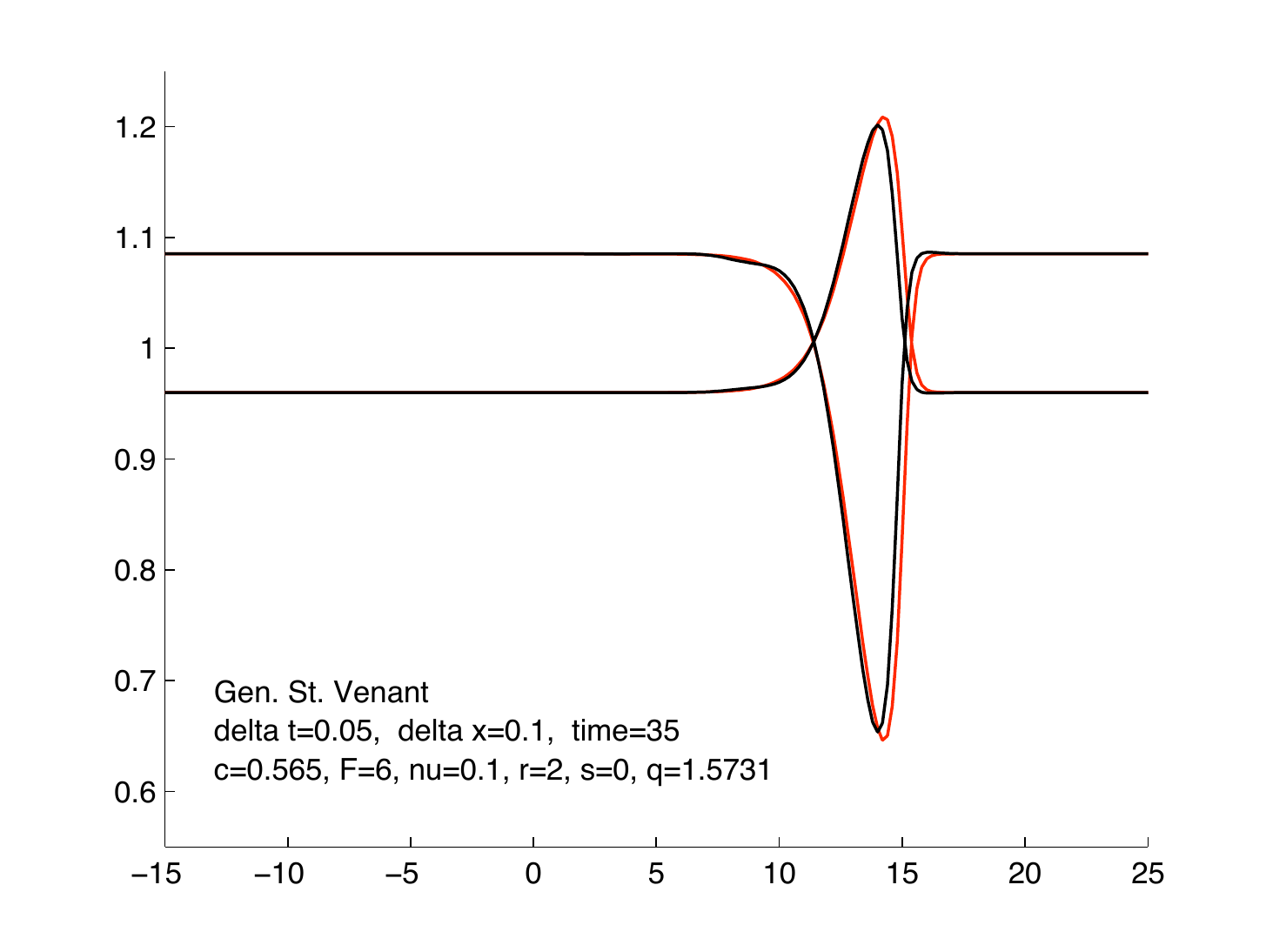}&\includegraphics[scale=0.25]{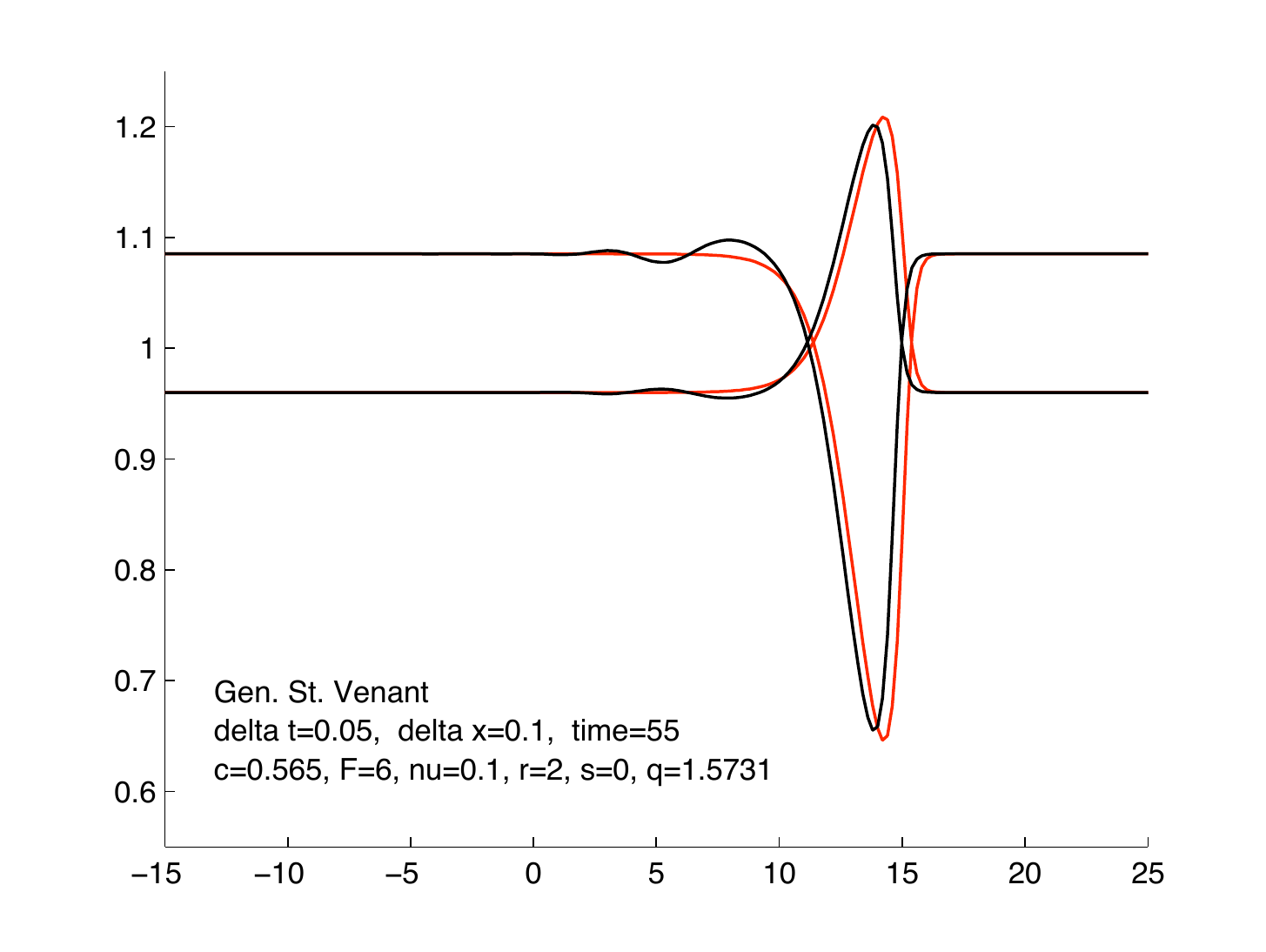}& \includegraphics[scale=.25]{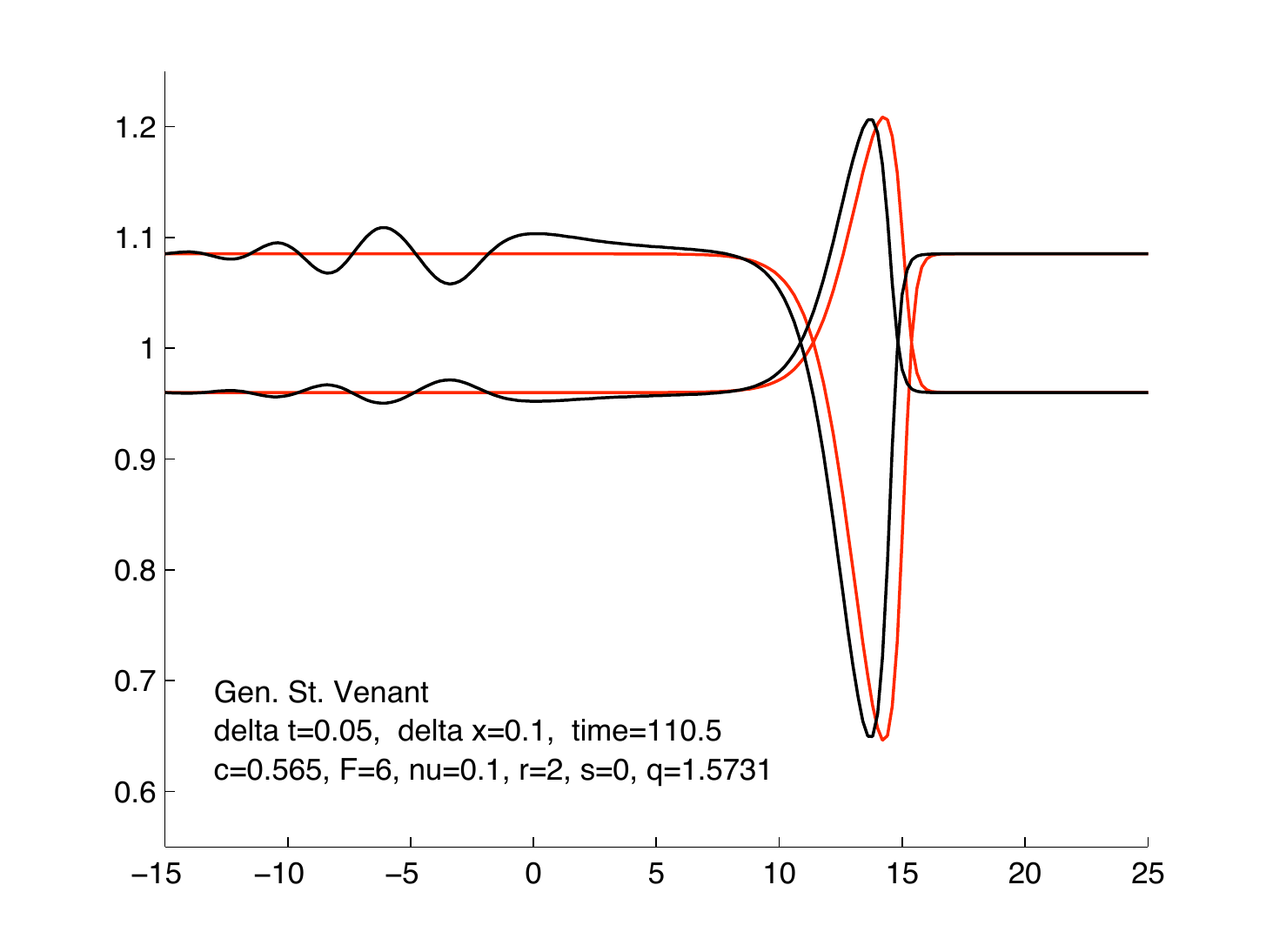}
\end{array}
$
\caption{Time evolution of square wave perturbation of a homoclinic orbit possessing stable point spectrum and unstable essential
spectrum. Here $u_-=0.96$, $q=u_-+c/u_-^2$, $\nu=0.1$, $r=2$, $s=0$, and $F=6$.  Notice that the perturbation decays
across the region of the profile where the gradient varies non-trivially and grows near the unstable limiting constant states
and is eventually convected to minus infinity.}
\label{f:evol}
\end{center}
\end{figure}

It follows that any spectrally stable periodic traveling wave of \eqref{e:stvcommon}, if one exists, must be of intermediate amplitude.
In particular, due to the complicated nature of the linearized operator and the fact that we can not rely on perturbation
techniques from a particular limit, our analysis now turns over to a numerical study.  Before continuing, however, we wish
to give a heuristic argument,
reconciling physically observed stability with this
analytically-demonstrated instability,
why one
still might
expect the existence of stable periodic waves of intermediate amplitude
despite the instability of the Hopf and homoclinic limiting states.
If one considers the stability of the limiting homoclinic profiles more
carefully, it can be found by rigorous Evans function calculations\footnote{
More precisely, numerical Evans function computations with rigorous
error bounds; see \cite{BJRZ}.}
that there exist homoclinic orbits having unstable essential spectrum, as predicted from the preceding discussion,
and \emph{stable} point spectrum.
As a result, we find that the associated homoclinic orbit stabilizes perturbations across
dynamic parts of the wave, i.e. where the gradient varies nontrivially, reflecting the stable point spectrum of the wave,
while the portion of the wave near the limiting constant endstate amplifies the perturbation, reflecting the unstable
essential spectrum; see Figure \ref{f:evol}.
Accordingly, one encounters an interesting ``metastability" mechanism where the stable point
spectrum induces a stabilizing effect
on
a closely spaced array of solitary waves, since the unstable constant endstates
would have little effect due to the ``closeness" of the array.  This leads one to a notion of the ``dynamic stability" of a solitary
wave, which is essentially the spectrum of an appropriately periodically extended version of the original homoclinic;
this issue is discussed in more detail in \cite{BJRZ}.
Heuristically, then, considering a closely spaced array
of solitary solutions as a periodic orbit we are led to the possibility of finding spectrally stable periodic waves away
from either the homoclinic or Hopf limits.  In the next section, we numerically verify this heuristic by presenting
numerical computations indicating
the existence of a stable periodic solution to the St. Venant equations \eqref{e:stvl}.
These numerical stability results are
formal in the sense that we do not present any error bounds or rigorous
high-frequency asymptotics precluding the existence of unstable spectrum sufficiently far from the origin; in a future
paper \cite{BJNRZ} we hope to present such bounds and hence make the formal numerical
arguments here rigorous.  
For the purposes of this article, however, our formal numerical investigation will suffice.


%
%
%

\subsection{Numerical study}\label{s:num}

\begin{figure}[htbp]

\begin{center}
$
\begin{array}{ccc}
  \includegraphics[scale=.25]{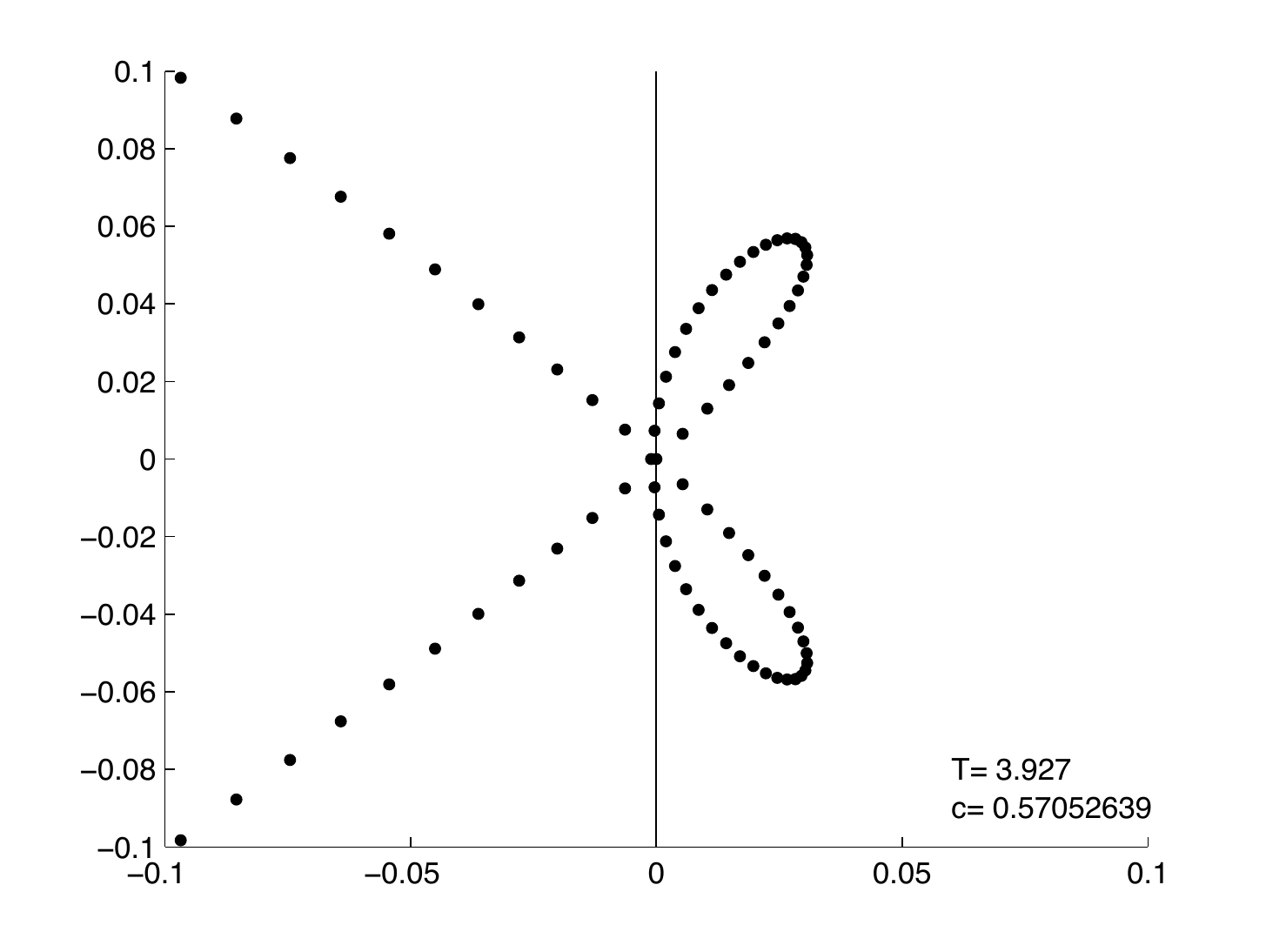}&\includegraphics[scale=0.25]{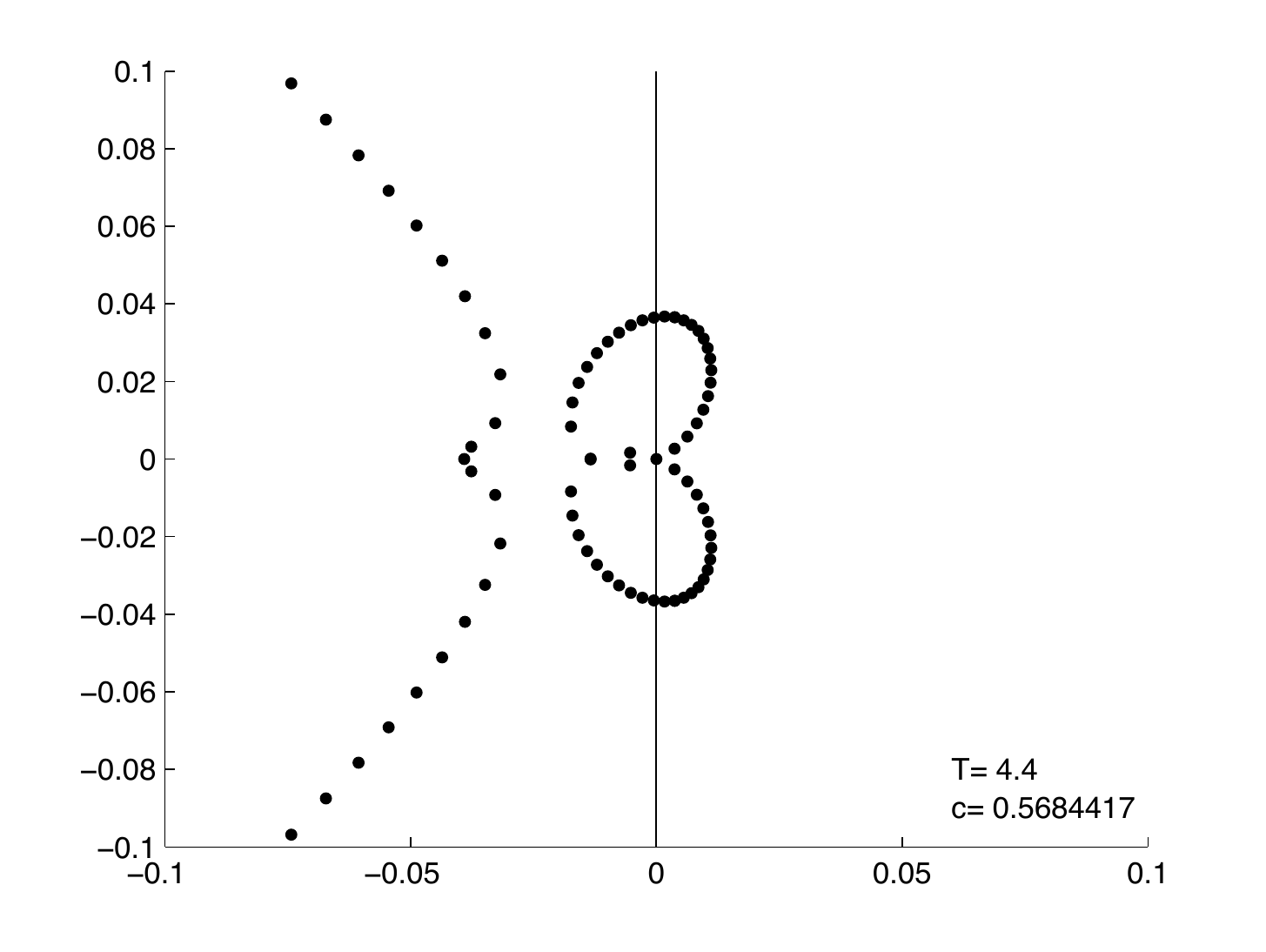}& \includegraphics[scale=.25]{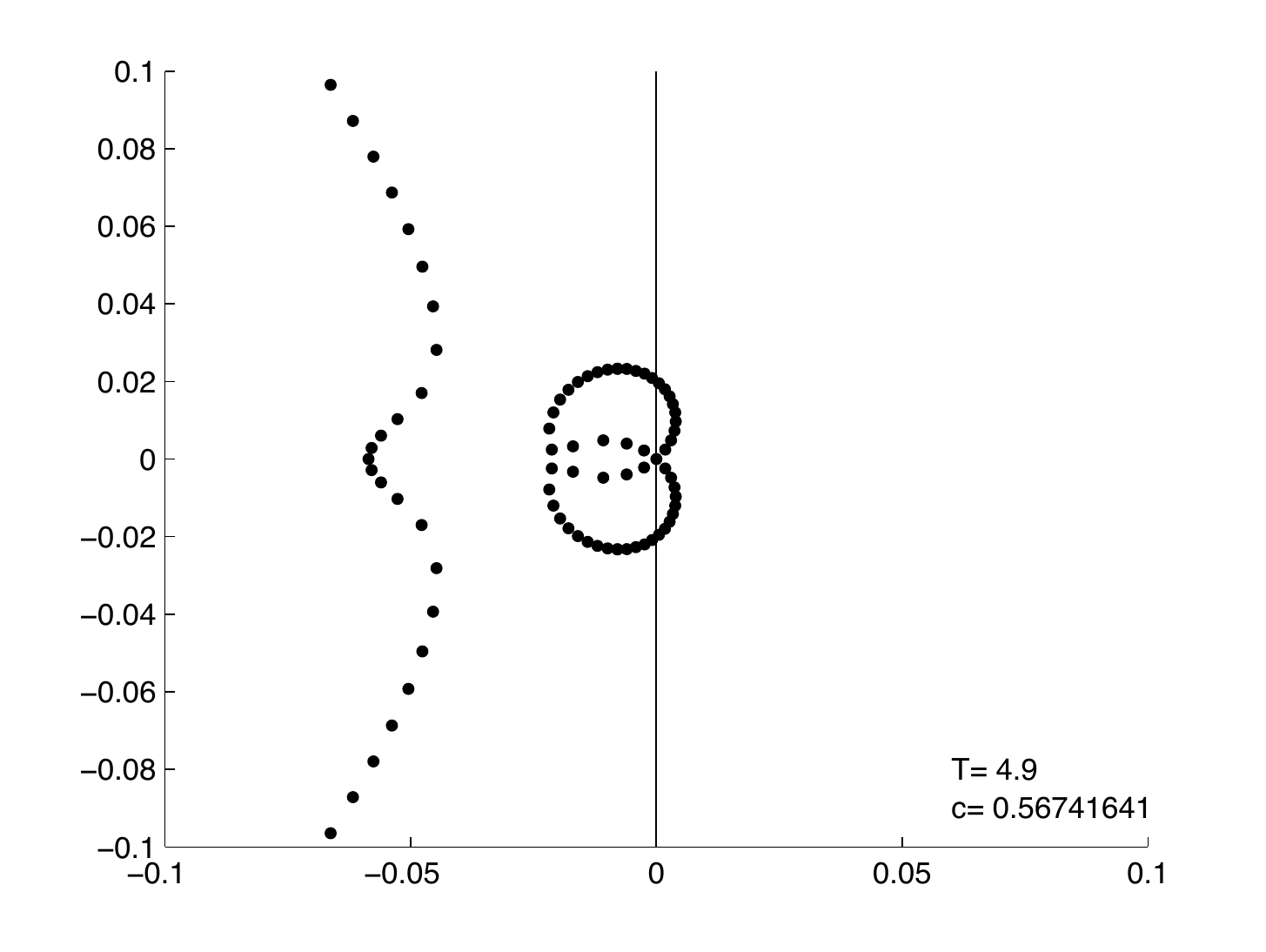}\\
  \includegraphics[scale=.25]{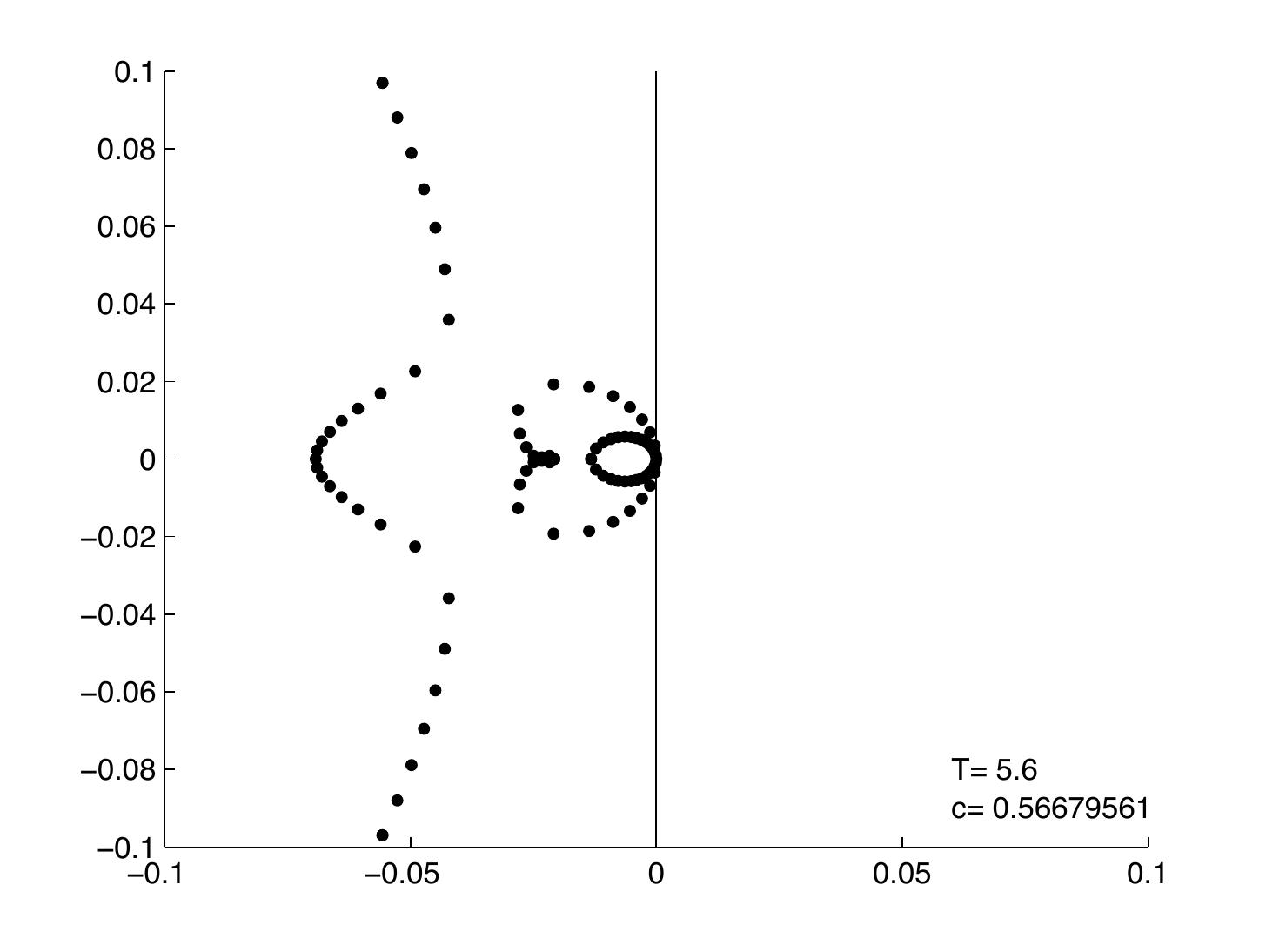}&\includegraphics[scale=0.25]{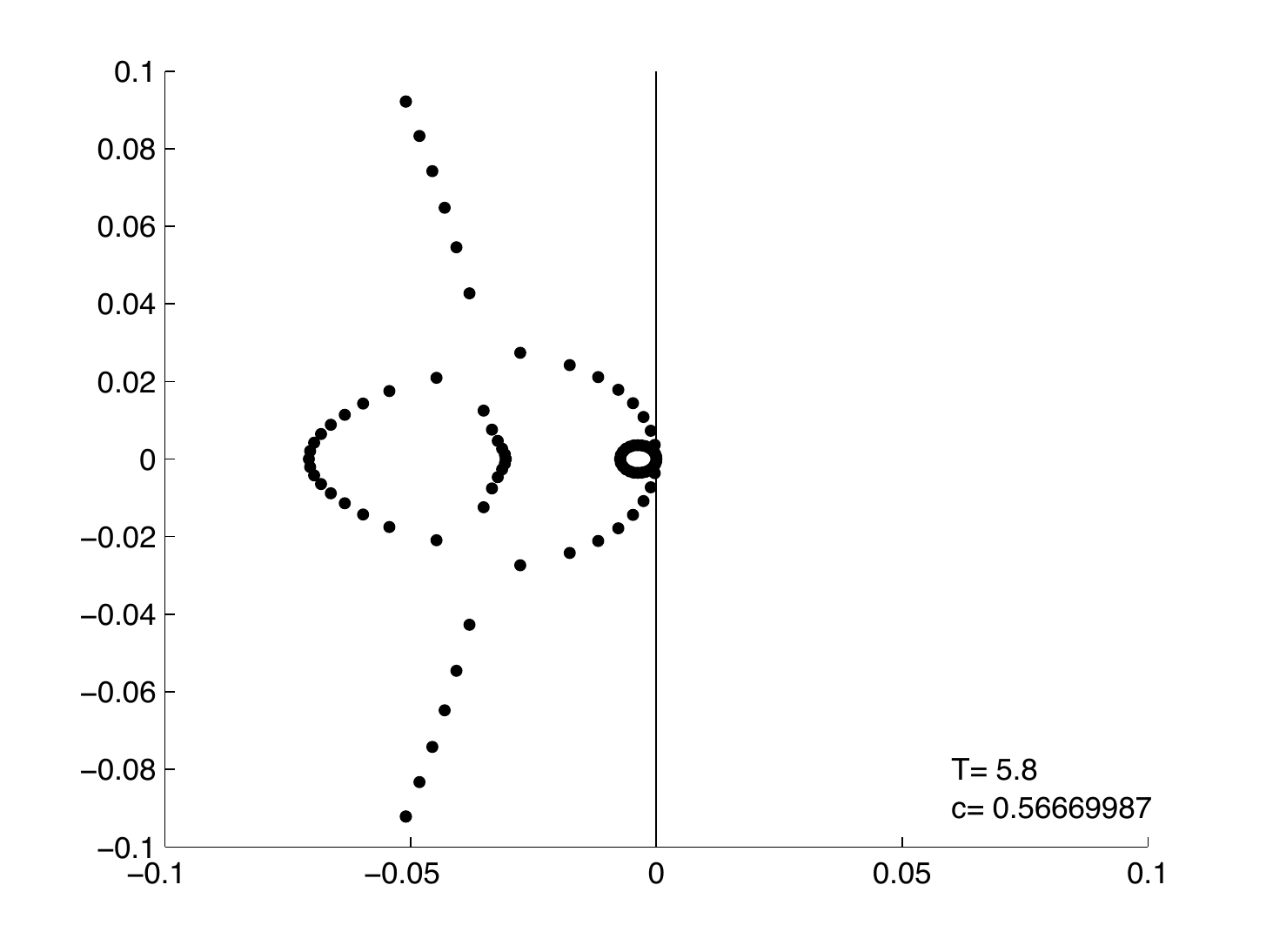}& \includegraphics[scale=.25]{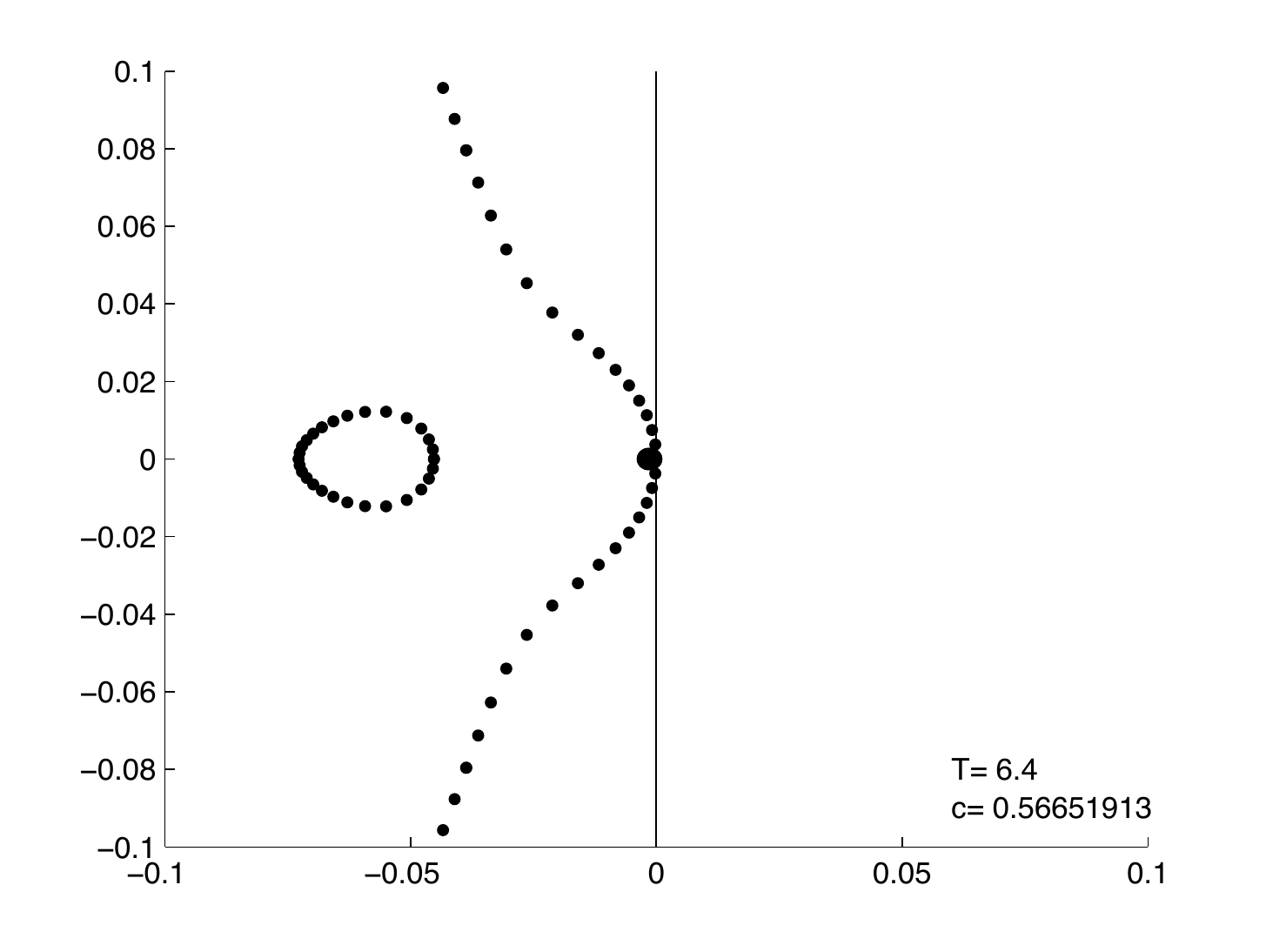}\\
  \includegraphics[scale=.25]{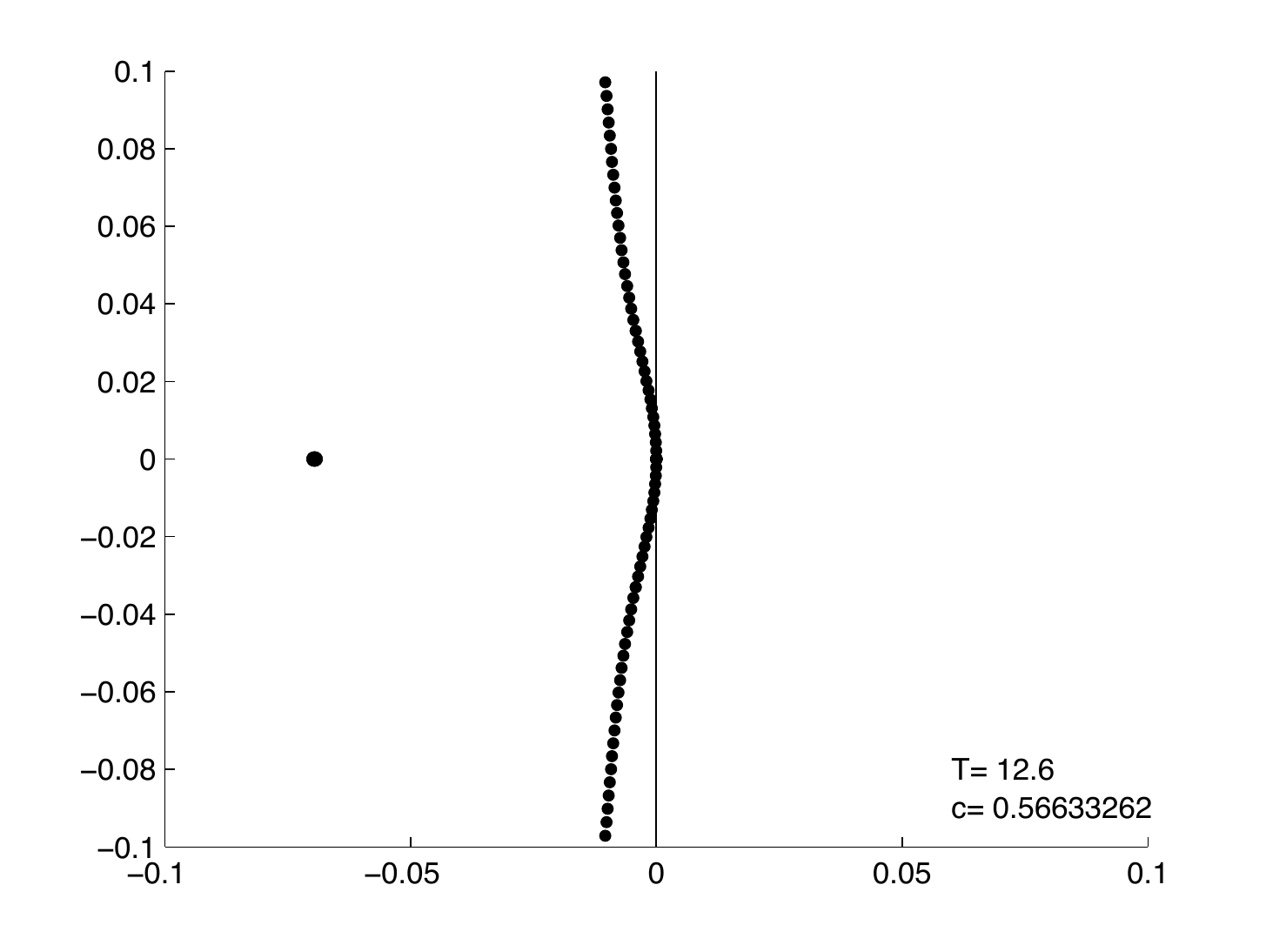}&\includegraphics[scale=0.25]{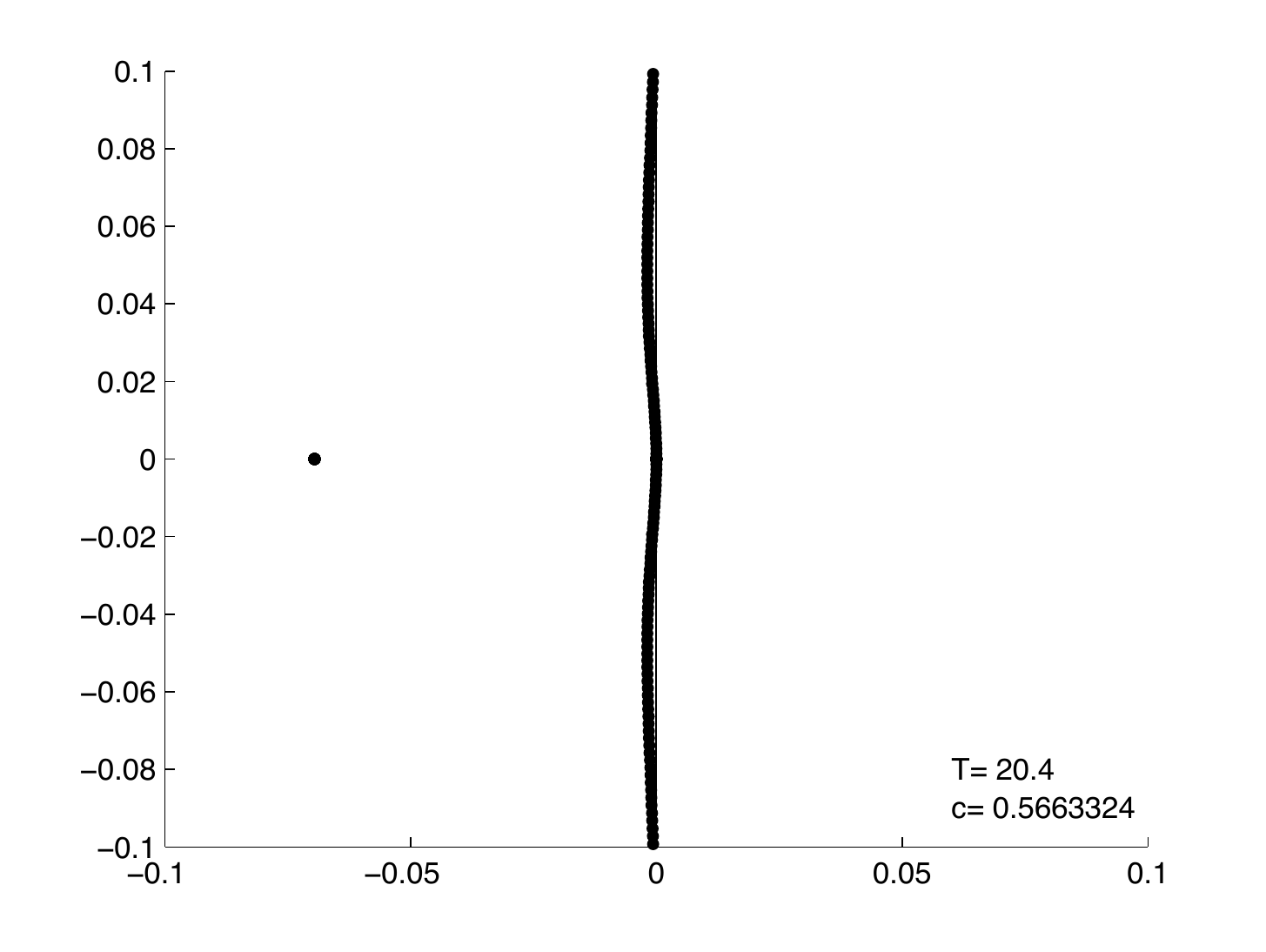}& \includegraphics[scale=.25]{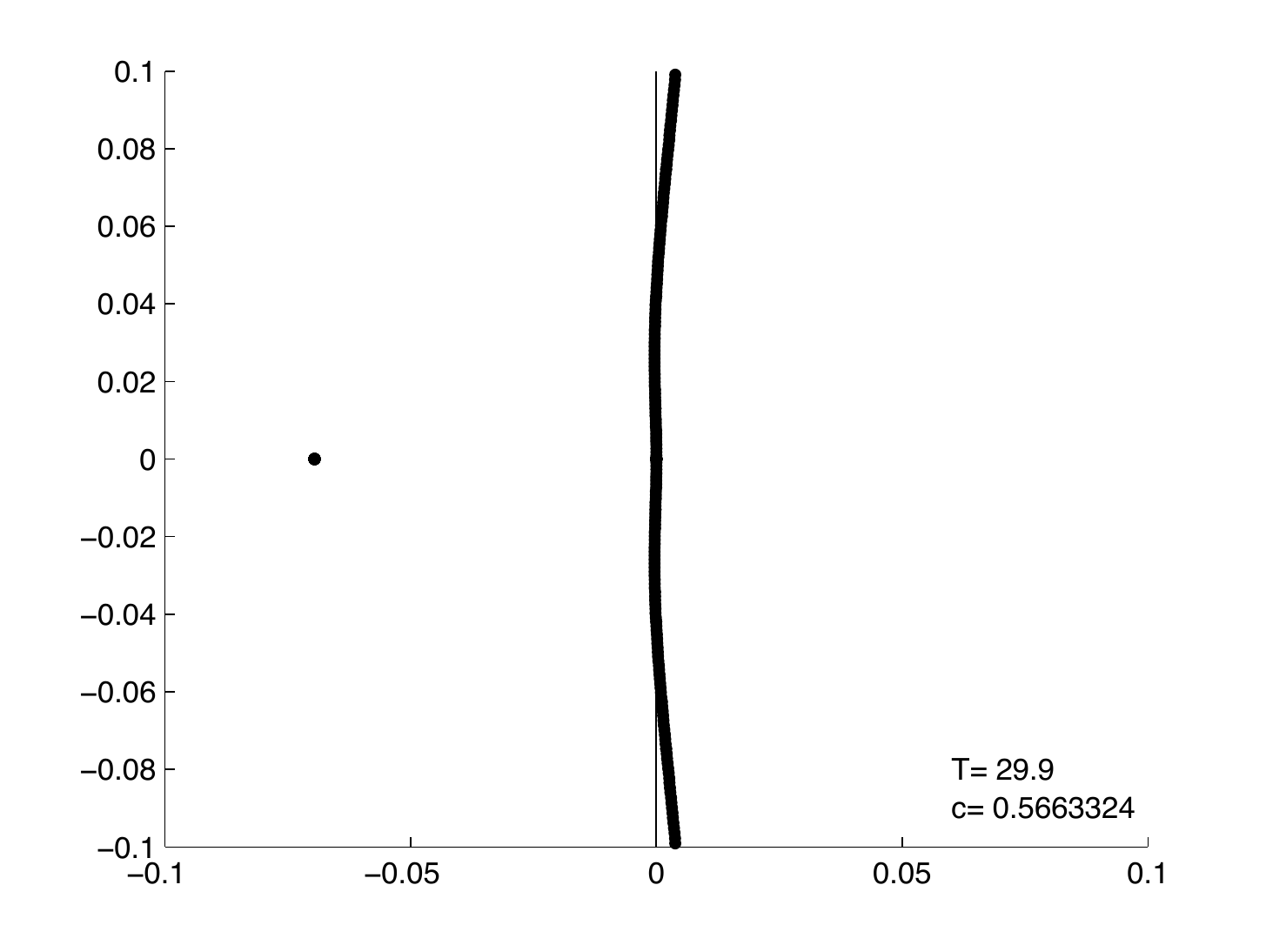}
\end{array}
$
\end{center}
\caption{Evolution of spectra as period, here denoted as $X=T$, and wave speed, $c$, vary. Here $u_-=0.96$, $q=u_-+c/u_-^2$, $\nu=0.1$, $r=2$, $s=0$, and $F=6$.
Starting in the top left picture and running from left to right and from top to bottom, we see the evolution of the spectra from the Hopf
bifurcation at $T\approx 3.9$ to a wave seemingly near the the homoclinic in phase space with period $T\approx 29.9$.}
\label{f:animate}
\end{figure}

In this section, we use the SpectrUW package, which relies on Fourier-Bloch decompositions and Galerkin  truncation to
numerically evaluate the spectrum of linear operators with periodic coefficients, to numerically compute
the spectrum of the periodic orbits depicted in Figure \ref{f:phase}.  As described in the previous section, we expect
the solutions near the Hopf and homoclinic cycles to have unstable essential spectrum.  Nevertheless, the metastability
of the limiting homoclinic profile suggests that waves of intermediate period may be spectrally stable, and hence
be nonlinearly stable by Theorem \ref{t:nl}.  In this investigation then, we animate the spectrum
as the period $X$ is increased, or equivalently the wave speed $c$ is decreased, from the Hopf period $X\approx 3.9$
to $X=29.9$, which corresponds to a periodic orbit seemingly very close to the homoclinic in phase space.  The results
of this study are shown in Figure \ref{f:animate} and indeed seem to indicate a region of stable periodic orbits.
Indeed, it seems that the unstable small-amplitude periodic waves eventually stabilize in a neighborhood of the origin as the period is increased
and then are later destabilized by the essential spectrum crossing the imaginary axis at non-trivial complex conjugate points as the period is increased
further.

\begin{figure}\begin{center}
$\begin{array}{cc}
\includegraphics[trim=0 0 0 0, scale=.3]{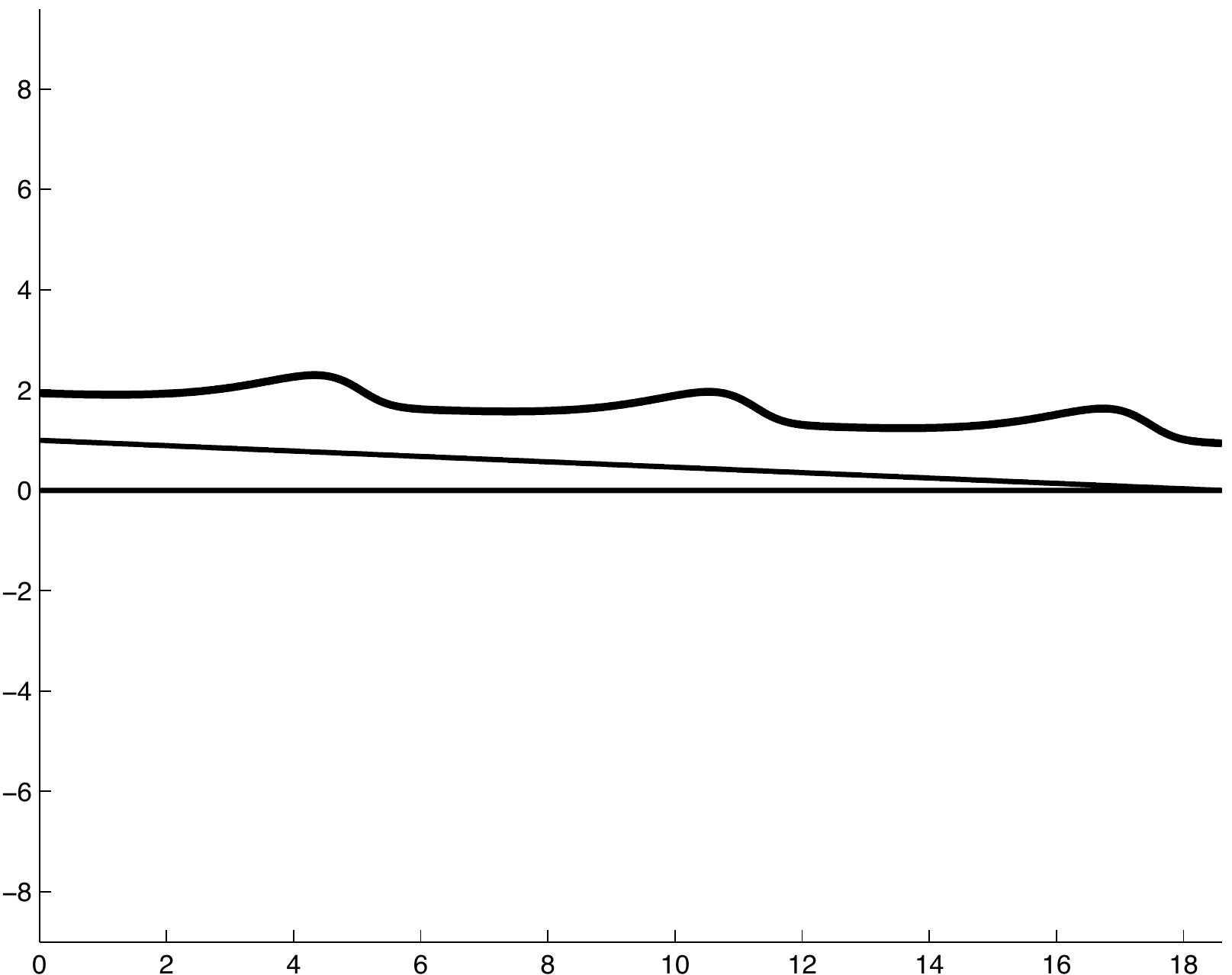}
\end{array}$
\caption{A numerically
stable periodic wave of the St. Venant equation \ref{e:stvcommon}, plotted
in the original physical coordinates ($h=\tau^{-1}$ vs. $x$).  This particular wave has period $X\approx 6.2$,
and corresponds to the bold orbit in Figure \ref{f:phase}(a) between the upper and lower stability boundaries,
and whose period is designated by the circle in Figure \ref{f:phase}(b) .  In particular,
notice that the corresponding wave speed, as
depicted in Figure \ref{f:phase}(b), is close to the limiting homoclinic speed.}
\label{f:stable}
\end{center}
\end{figure}

We see then that there seems to be a regime of stability in which periodic orbits with particular intermediate periods
are spectrally stable solutions of \eqref{e:stvcommon}.  A spatial plot in the original physical
coordinates ($h=\tau^{-1}$ vs. $x$) of a periodic roll-wave in this stable regime is depicted in
Figure \ref{f:stable}.  In \cite{BJNRZ}, high frequency asymptotics have been obtained
which make this numerical evidence rigorous by proving that any spectral instability must occur within a specified compact
region of the complex plane.  Furthermore, for the seemingly stable spectra depicted in Figure \ref{f:animate} it is
verified in \cite{BJNRZ} through the development of rigorous error bounds that the corresponding waves are indeed spectrally stable as solutions of \eqref{e:stvcommon}.

Finally, we make some remarks concerning the various instabilities present in Figure \ref{f:animate} and their relation to
the hyperbolicity of the associated Whitham averaged system.  To begin, we recall that by the recent work \cite{NR}
hyperbolicity of this Whitham system is necessary for spectral stability; see Theorem \ref{t:nr}.  The lack of sufficiency
in this theorem is associated with the fact that it is only a \emph{first order} verification.  Hence, Theorem \ref{t:nr}
essentially states that hyperbolicity of the Whitham averaged system is equivalent to the spectrum of the associated
linearized spectral problem being tangent to the imaginary axis at the origin, which is clearly a necessary condition
for stability but not sufficient\footnote{This is in contrast to the dispersive Hamiltonian case, such as the generalized
Korteweg-de Vries or nonlinear Schrodinger equations, where the Hamiltonian structure of the linearized operator implies
the stability spectrum is invariant with respect to reflections across the imaginary axis.  In that case, generically
it can be shown that hyperbolicity of the Whitham equation is \emph{equivalent} with spectral stability of the underlying
wave \emph{in a neighborhood of the origin}.  See, for example, \cite{JZ1}.}.  As an example, notice that the first
three spectral plots, ordered from left to right and top to bottom, in Figure \ref{f:animate} are spectrally unstable
in a neighborhood of the origin due to lack of hyperbolicity of the associated Whitham averaged system.  The remaining six
spectral plots are seemingly associated with hyperbolic Whitham averaged systems, but we see instability arising for
sufficiently large period due to an essential instability occurring away from the origin.
Thus, as expected, hyperbolicity of this Whitham equation is only a \emph{local} condition for spectral stability,
in the sense it only detects instabilities in a neighborhood of the origin.
In particular, the lower stability boundary, occurring 
with period within $0.1$ of $X=5.3$ is marked by the hyperbolic
Whitham criterion, while the upper stability boundary, occuring around $X=20.6$ is not.

Notice, however, that in the final spectral plot the wave seems to stabilize in a neighborhood of the origin.  While this seems to be a general phenonamon for periodic
waves where the period is not ``too" large, tentative numerical experiments indicate that for periodic waves with
with very large periods the spectrum seems to destabilize in a neighborhood of the origin; in particular, the spectrum
eventually seems to resemble that of Figure \ref{f:cst}(b) for the limiting homoclinic.  This seems to suggest that,
although the Whitham averaged system is hyperbolic, the wave is spectrally \emph{unstable} in a neighborhood
of the origin.  Readers should be warned, however, that a stabilization effect near the origin may still occur, but we may
not be able to see it due to a low-resolution of the spectral plot.  Furthermore, it seems to be quite difficult
to numerically generate periodic orbits with very large period and hence we had to resort to periodically extending
a homoclinic orbit for these numerics.  Nevertheless, these experiments seem to suggest that it may be possible
for a periodic traveling wave solution of the St. Venant equations \eqref{e:stvcommon} with sufficiently large
period to have a hyperbolic Whitham averaged system but be spectrally unstable to long-wavelength perturbations.

Continuing, we should note that using a Bloch-wave expansion, two of the authors of the present paper have recently
been able to rigorously validate the \emph{second order} Whitham expansion \cite{NR}, showing that this second order
Whitham system determines the convexity of the spectrum near the origin: a clearly more refined feature
than the first order verification provided by hyperbolicity.  Thus, it may be possible to numerically verify
the existence of periodic waves where the spectrum leaves the origin at second order and moves into the unstable
half plane by analyzing the associated second order Whitham system.  This  system is however considerably
more complicated than its first order counterpart discussed in Section \ref{s:whitham} and it is not immediately
clear how to ``compute" the necessary information from this system for a given numerically generated solution.

\section{Conclusions \& Discussion}

In this note, we have considered both analytical and numerical aspects of the stability of periodic roll-wave solutions
of the generalized St. Venant equations.  In particular, we reviewed known results concerning the nonlinear stability
of such solutions and then proceeded to numerically investigate the necessary spectral stability assumptions in the
nonlinear stability theorem.  To this end, we utilized the SpectrUW package developed at the University of Washington
and, formally, we made the case for the existence of a spectrally, and hence nonlinearly, stable periodic
traveling wave solution of the governing St. Venant equation.  This stands in contrast to the fact that periodic
solutions near either the Hopf equilibrium solution or the bounding homoclinic solutions are unstable.
By briefly outlining the heuristic of the ``dynamic stability" of the bounding homoclinic wave, however, we were
able to give a (possibly general) explanation of how an equation with unstable solitary waves can admit
stable periodic waves solutions.

This concept of ``metastability" of solitary waves has been considered also by Pego, Schneider, and Uecker \cite{PSU} in
the context of the related fourth-order diffusive Kuramoto--Sivashinsky model
\[
u_t+\partial_x^4u+\partial_x^2u+\frac{\partial_x u^2}{2}=0,
\]
which has been proposed as an alternative model for thin film flow down an inclined ramp.  Therein, the authors
analyze the time-asymptotic behavior of solutions of this equation, and conclude that they are dominated by trains
of solitary pulses.  As such, the mechanism of stable ``dynamic spectrum" seems to provide a partial
answer for how a train of solitary pulses can stabilize the convective instabilities shed from their neighbors.

In general, however, it seems that the mechanism
by which an equation with an unstable solitary wave can admit
stable periodic waves is not completely clear,
although we suspect that it is closely tied with the notion of the
``dynamic stability" of the limiting homoclinic profile.
In particular, to prove a general theorem describing this
mechanism at a quantitative level
seems outside the realm of current methods and
remains an interesting open problem.

Finally, we wish to emphasize again
that the numerical evidence for stability presented
in Section \ref{s:num} are formal in the sense that they are lacking the error estimates and high-frequency
asymptotics/energy estimates necessary to preclude the existence of unstable spectrum outside the window of our computation.
In a future paper \cite{BJNRZ}, we will carry out these details and provide rigorous numerics which indicate the
existence of a spectrally stable periodic traveling wave solution of the St. Venant equation \eqref{e:stvcommon}.

\medskip
{\bf Acknowledgement}: Thanks to
Bernard Deconink for his generous help in guiding us in the
use of the SpectrUW package developed by him and collaborators.
The numerical Evans function computations performed
in this paper were carried out
using the STABLAB package developed by Jeffrey Humpherys with
help of the first and last authors.

\end{document}